\documentclass[12pt]{article}
\usepackage{graphicx}

\def\F{ {\cal F} }
\def\pd#1{ \partial_{#1}}
\def\QED{ Q.E.D. }
\newtheorem{theorem}{Theorem}
\newtheorem{proposition}{Proposition}
\usepackage[usenames,dvipsnames]{color}

\newcommand{\nobukiNote}[1]{ }

\usepackage{amsfonts, amsmath}
\usepackage{amssymb}
\newtheorem{definition}{Definition}
\newtheorem{lemma}{Lemma}
\newtheorem{corollary}{Corollary}
\newtheorem{remark}{Remark}
\newtheorem{example}{Example}
\newcommand{\CC}{{\mathbb C}}
\newcommand{\RR}{{\mathbb R}}
\newcommand{\QQ}{{\mathbb Q}}
\newcommand{\PP}{{\mathbb P}}
\newcommand{\ZZ}{{\mathbb Z}}
\newcommand{\NN}{{\mathbb N}}
\newcommand{\cF}{{\mathcal F}}
\newcommand{\cH}{{\mathcal H}}

\newcommand{\cD}{{\mathcal D}}
\newcommand{\cO}{{\mathcal O}}
\newcommand{\cM}{{\mathcal M}}

\newcommand{\cQ}{{\mathcal Q}}

\newcommand{\gr}{{\rm gr}}
\newcommand{\wA}{{{\widetilde A}(w)}}
\newcommand{\wa}{{\widetilde{a}}}
\newcommand{\inn}{{\rm in}}

\usepackage{datetime}
\usepackage{times}

\title{Irregular Modified $A$-Hypergeometric Systems}
\author{Francisco-Jes\'us Castro-Jim\'enez\thanks{First two authors are partially supported by MTM2010-19336 and
FEDER and Junta de Andaluc\'{\i}a FQM5849, FQM333. FJCJ is also partially supported by S-13025-JSPS (Japan); MCFF is also partially
supported by Max Planck Institute f\"ur Mathematik (Bonn). Third author is partially supported by JSPS grants-in-aid No. 21740098 and No. S-24224001.},
Mar\'{\i}a-Cruz Fern\'andez-Fern\'andez$^*$, \\
Tatsuya Koike and Nobuki Takayama}
\date{}
\begin{document}
\maketitle
\begin{center}
\today
\end{center}

\begin{abstract}
A modified $A$-hypergeometric system is a system of differential equations
for the function $f(t^w \cdot x)$ where $f(y)$ is a solution of an
$A$-hypergeome\-tric system in $n$ variables and $w$ is an $n$ dimensional
integer vector, which is called the weight vector. We study the irregularity of modified systems
by adapting to this case the notion of umbrella introduced by M.
Schulze and U. Walther. Especially, we study slopes and Gevrey
series solutions. We develop  some applications of this study. Under some conditions
we give Laplace integral representations of divergent series
solutions of the modified system and we show that certain Gevrey series
solutions of the original $A$-hypergeometric system along
coordinate varieties are Gevrey asymptotic expansions of holomorphic solutions of the $A$-hypergeometric system.
\end{abstract}

\tableofcontents

\section{Introduction}

$A$-Hypergeometric systems (or GKZ-systems or simply hypergeometric systems) are systems of linear partial differential equations on the complex affine space $\CC^n$.
Although they were already considered in works by J.~Hrabowski \cite{Hra85} and by I.~M.~Gelfand, M.~I.~Graev and A.~V.~Zelevinsky \cite{ggz1987},
the systematic study of hypergeometric systems started with the paper by I.~M.~Gel'fand, M.~M.~Kapranov and A.~V.~Zelevinsky
\cite{gkz1989}. Each of these systems, denoted by $H_A(\beta)$, is determined by a pair $(A,\beta)$ where $A=(a_{ij})$ is an integer $d\times n$ matrix of rank $d$ and $\beta\in \CC^d$ is a parameter vector.
An $A$-hypergeometric system $H_A (\beta)$ is defined by the $d$ Euler operators $E_i-\beta_i:=\sum_{j=1}^n a_{ij}y_j\frac{\partial}{\partial y_j}-\beta_i$ for $i=1,\ldots,d$ and the toric operators $\partial^u-\partial^v$ associated with each pair $(u,v)\in \NN^n\times \NN^n$ such that $Au=Av$. Here $\partial^u$ stands for the monomial differential operator $\partial_1^{u_1}\cdots\partial_n^{u_n}$ and $\partial_i = \frac{\partial}{\partial y_i}$.

For generic parameters $\beta \in \CC^d$, the holomorphic solutions of $H_A(\beta)$ at nonsingular points can be described by using the so-called $\Gamma$--hypergeometric series (see \cite{ggz1987}, \cite{gkz1989}; see also \cite{SST} and Subsection \ref{Gevrey-solutions-at-infinity}).
We are interested in divergent $\Gamma$--hypergeometric series solutions. To study them we use the notion of slopes defined in the general setting by Y. Laurent \cite{Laurent}.
If the matrix $A$ is pointed (which means that the column vectors of ${A}$  lie  in an open half-space with boundary passing
through the origin in $\RR^{d}$), M. Schulze and U. Walther \cite{SW} have described the slopes of $H_A(\beta)$ with respect to coordinates varieties, generalizing previous work in \cite{castro-takayama}, \cite{hartillo1} and \cite{hartillo2}. These slopes are closely related to the irregularity  of the system \cite{Laurent-Mebkhout-pa-pa}
and the existence of non convergent Gevrey series solutions of this system. Gevrey series solutions of $H_A(\beta)$ were studied by the second author in \cite{Fer} (see also \cite{fer-cas-1} and \cite{fer-cas-2} where the authors treat particular cases).

Modified hypergeometric systems were introduced by the fourth author \cite{takayama2009} in order to study solutions of hypergeometric systems along a curve $y(t) =(c_1 t^{w_1}, \dots , c_n t^{w_n})$ for $w\in \ZZ^n$, $c_i \in \CC$. Each of them
is determined by  a tuple $(A,w,\beta,\alpha)$ where $A$ and $\beta$ are as before, $w=(w_1, \ldots, w_n) \in \ZZ^n$ and $\alpha \in \CC$.
We denote by ${\widetilde A}(w)$ (or simply $\widetilde A$) the matrix
$${\widetilde A}(w) = {\widetilde A}=
\left( \begin{array}{cccc}
            a_{11} & \cdots & a_{1n} & 0 \\
                   & \cdots &        & 0 \\
            a_{dn} & \cdots & a_{dn} & 0 \\
            w_1    & \cdots & w_n    & 1 \\
\end{array} \right).
$$
Throughout this paper, we do not always assume that $A$ is pointed,
but we assume that ${\widetilde A}$ is.
Note that when $A$ is pointed, then ${\widetilde A}$ also is.

\begin{definition} \rm (\cite{takayama2009}) \label{def:mgkz}
We call the following system of differential equations $H_{A,w,\alpha}(\beta)$
a {\it modified  $A$-hypergeometric system}:
\begin{eqnarray}
  \left( \sum_{j=1}^n a_{ij}  x_j \partial_{j} - \beta_i \right) \bullet f &=& 0,
   \qquad(i = 1, \ldots, d)  \label{eq:mgkz1} \\
  \left( \sum_{j=1}^n w_j x_j \partial_{j} - t \partial_t -\alpha\right) \bullet f &=& 0,
   \label{eq:mgkz2} \\
  \left( \prod_{i=1}^n \partial_{i}^{u_{i}} t^{u_{n+1}}
       - \prod_{j=1}^n \partial_{j}^{v_{j}} t^{v_{n+1}}
  \right) \bullet f &=& 0 \label{eq:mtoric} \\
    \nonumber
\end{eqnarray}
with  $u, v \in {\NN}^{n+1}$ running over all $u, v$ such that ${\widetilde A}u = {\widetilde A}v$. Here we denote $\frac{\partial}{\partial x_i}=\partial_i$.
\end{definition}
The modified  system is defined on the space $X=\CC^{n+1}$ with coordinates $(x,t)=(x_1,\ldots,x_n,t)$.
Let $D$ (or $D_{n+1}$) be the Weyl algebra in $(x,t)$. The left ideal
in $D$ generated by the operators in (\ref{eq:mgkz1}),
(\ref{eq:mgkz2}), and  (\ref{eq:mtoric}) is also denoted by
$H_{A,w,\alpha}(\beta)$ if no confusion arises. The left $D$-module
$D/H_{A,w,\alpha}(\beta)$ is denoted by $M_{A,w,\alpha}(\beta)$.
By \cite{takayama2009}  the
$D$-module $M_{A,w,\alpha}(\beta)$ is holonomic for any
$A,\beta,w,\alpha$.

The system $H_A(\beta)$ is a summand of the modified system on the space $t\not= 0$.
More precisely, denote $Y=\CC^{n+1}$ with coordinates $(y,s)$ and consider the map
\begin{equation} \label{eq:gr-deformation}
\varphi
: Y^*:= \CC^n\times \CC^* \subset Y
  \longrightarrow  X^*:= \CC^n\times \CC^* \subset X
\end{equation}
defined by $\varphi
(y_1,\ldots,y_n,s)
=(s^{-w_1}y_1,\ldots,s^{-w_n}y_n,s)$. The pullback image of the ideal $H_{A,w,\alpha}(\beta)$
by $\varphi$ equals
the ideal of differential operators on
$Y^*$ generated by $H_A(\beta)$ and
$s\partial_s+\alpha$.
Notice that this last ideal is nothing but the hypergeometric ideal
associated with the matrix
$\widetilde{A} ({\bf 0})$
and the parameter vector $(\beta,-\alpha)\in \CC^{d+1}$.
As usual we denote this ideal by $H_{\widetilde{A}({\bf 0})}(\beta , -\alpha )$.

Let us consider the local analytic situation. Let ${\cD}$ be the sheaf of holomorphic differential operators on
$X$ and $\cM_{A,w,\alpha}(\beta)$ the quotient sheaf ${\cD}/{\cD}H_{A,w,\alpha}(\beta)$. By the previous observation,
the hypergeometric ${\cD}$-module $\cM_{\widetilde{A}({\bf
0})}(\beta ,-\alpha )={\cD}/{\cD} H_{\widetilde{A} ({\bf
0})}(\beta ,-\alpha )$ is the extension to $Y$ of the pullback
module $\varphi^*(\cM_{A,w,\alpha}(\beta))$ considered as a
$\cD_{Y^*}$--module on the first space $Y^*=\CC^n\times \CC^*$.
Since $\varphi$ is a biholomorphic map both $\cD$-modules
$\cM_{\widetilde{A}({\bf 0})}(\beta , -\alpha )_{|Y^*}$ and
$\cM_{A,w,\alpha}(\beta)_{|X^*}$ are isomorphic.

We can describe solutions of the original hypergeometric system $H_A(\beta)$ associated to
the weight vector $w\in \ZZ^n$ via solutions of the modified system.
Series solutions of $H_A(\beta)$ have
been studied in \cite{ggz1987} and \cite{gkz1989} where the authors constructed
convergent series solutions associated to the regular triangulation
induced by a generic weight vector $w$. The
construction is generalized as follows \cite{SST}: Assume
that $\beta\in \CC^d$ is very generic (this condition is essential in the construction). Suppose that the initial ideal ${\rm in}_{(-w,w)}(H_A(\beta))$ has a solution of the form $y^\rho$,
$\rho \in \CC^n$. Then the monomial $y^\rho$ can be extended
to a formal series solution $\phi(y) = y^\rho + \cdots$ of
$H_A(\beta)$. We call the series $\phi(y)$ a series solution of $H_A(\beta)$
associated to the weight vector $w$. The series is divergent in
general. We are interested in giving an explicit expression of a
solution of $H_A(\beta)$ whose  asymptotic expansion is $\phi(y)$. A
standard method, in the theory of ordinary differential equations,
to construct such an expression is the Laplace integral representation
and the Borel transformation of
divergent series. This method has been successful in the
study of global analytic properties of solutions of ordinary
differential equations, see, e.g., the book by W. Balser \cite{B} and
the references therein. Then it is a natural problem to construct a Laplace
integral representation corresponding to the divergent series solution
$\phi(y)$ associated to the weight vector $w$. Our modified system,
which is a system of differential equations for
$\phi(t^{w_1}x_1, \ldots,  t^{w_n}x_n)$, is used to give an answer
to this problem.

A key ingredient of our study is the fact that
the modified hypergeometric system
$H_{A,w,\alpha}(\beta)$ is transformed into
the hypergeometric system $H_{\widetilde{A}(w)}(\beta , \alpha -1 )$ associated with the matrix $\wA$ and parameter $(\beta ,\alpha -1 )$,
by the formal inverse Fourier transform
$t \mapsto \pd{t}$, $ \pd{t} \mapsto -t$,
which enables us to study the modified system using the theory of hypergeometric systems
(see Subsection \ref{section-Fourier}).

For example,
for $A=(1,2)$ and $\beta\in \CC$
the
system $H_A(\beta)$ is generated
by $x_1 \pd{1} + 2 x_2 \pd{2}-\beta {\mbox{ and }} \pd{1}^2-\pd{2}$. The
modified system $H_{A,w,\alpha}(\beta)$ for $w=(-1,-1)$, $\alpha \in \CC$ is generated by the
three  operators $x_1 \pd{1} + 2 x_2 \pd{2}-\beta, \,\, -x_1 \pd{1}
-x_2 \pd{2}-t \pd{t}-\alpha, \,\, \pd{1}^2 t -\pd{2}$, and the inverse Fourier
transformation of $H_{A,w,\alpha}(\beta)$ is generated by $x_1 \pd{1} + 2
x_2 \pd{2}-\beta, \,\, -x_1 \pd{1} -x_2 \pd{2}+t \pd{t}-\alpha + 1, \,\, \pd{1}^2
\pd{t} -\pd{2}$ which is equal to $H_{\widetilde{A}(w)}(\beta ,
\alpha-1 )$.

We study the behavior of solutions of
$H_{A,w,\alpha}(\beta)$ near the hyperplane $t=0$ in the space
$X$ and we will give Laplace integral representations of its solutions.

The structure of the paper is as follows: in Section \ref{generalities-slopes} we recall Y. Laurent's
definition of the slopes of a finitely generated $D$-module with
respect to a hypersurface.

In Section \ref{section-irregularity-GKZ} we recall the use of
umbrellas for the description of the slopes of a hypergeometric
system given by Schulze and Walther \cite{SW}, then we summarize the
construction of the Gevrey solutions given in \cite{Fer} and extend
some of these results to the case of the Gevrey solutions at
infinity.

In Section \ref{section:irregularity-MGKZ} we prove that the formal
inverse Fourier transform, with respect to $T$,  of a modified
hypergeometric system is an $A$-hypergeometric system and we
use this fact to describe the slopes of the former by using the
umbrella of the latter. We construct, associated with any slope of the
modified system $\cM_{A,w,\alpha}(\beta)$, a basis of its Gevrey solutions, modulo convergent
power series, when the parameters $\beta$ and $\alpha$ are very
generic, see Theorem \ref{formal-modulo-convergent}.
Moreover, for $\beta\in \CC^d$ very generic we construct a
basis of formal series solutions of the modified system for any
$\alpha\in \CC$ and $w\in \ZZ^n$, see Theorem
\ref{theorem-dimension-1}. If in addition $w$ is generic this basis is reduced to a single element. Later, in Sections \ref{section:Borel}
and \ref{section:borel_transform_revisited}, we prove under some assumptions that
this solution is an asymptotic expansion of a Laplace integral
representation of a solution by the Borel summation method
(Theorem \ref{borel-thm} and Section
\ref{section:borel_transform_revisited}). As an application, we give
an asymptotic error evaluation of finite sums of formal series
solutions of original $A$-hypergeometric systems with irregular
singularities studied in, e.g., \cite{Fer} and \cite{gkz1989} (see
Theorem \ref{borel-thm}, the inequality
(\ref{eq:asymptotic_expansion}), and Section
\ref{section:borel_transform_revisited}). Example
\ref{example12} illustrates the application for the simplest $A$.

Several integral representations have been studied for solutions of
regular holonomic hypergeometric systems (see, e.g.,
\cite{gkz-advance}, \cite{beukers}, \cite{bfp}). They play a
prominent role in the study of $A$-hypergeometric functions.
However, there have been few studies of integral representations for
solutions of irregular $A$-hypergeometric systems. We would like to
point out that integral representations of holomorphic solutions of
irregular $A$-hypergeometric systems have been recently given by A.
Esterov and K. Takeuchi \cite{esterov-takeuchi} by using the so
called rapid decay homology cycles. Our Laplace integral
representation, which we propose in this paper for giving an
analytic meaning to divergent series solutions, is different from
their representation: the integrand of our representation is an
$A$-hypergeometric function associated to a homogenized
configuration of $A$ (Section \ref{section:Borel}).

In Section \ref{section:borel_transform_revisited} we illustrate how, under some conditions, the study of the irregularity of
$\cM_{A,w,\alpha}(\beta)$ along $T$ gives
an analytic meaning to the Gevrey series solutions of $\cM_A(\beta)$, along coordinate varieties, constructed in
\cite{Fer}. More precisely we prove (see Proposition
\ref{bounded-growth}) that  they are asymptotic expansions of
certain holomorphic solutions of $\cM_A(\beta)$.

An interplay of algebra (slopes and formal series) and analysis
(Borel summation method) is a main point of this paper.
In order to make a comprehensible presentation to readers from
several disciplines, we often review some well-known facts to
experts. We hope that our style is successful.

Acknowledgements: We wish to thank J. Gonz\'{a}lez-Meneses, M. Granger and D. Mond for their help, suggestions and comments. We are very grateful to an anonymous referee, whose thoughtful suggestions have improved this article.

\section{Generalities on slopes} \label{generalities-slopes}

Recall that $D=D_{n+1}$ is the Weyl algebra
${\CC}\langle x_1, \ldots, x_n, t,\pd{1}, \ldots, \pd{n}, \pd{t}
\rangle$. In this section we write $x=(x_1,\ldots,x_{n+1})$, $\partial =(\partial_1,\ldots,\partial_{n+1})$. The variable $t$ is also denoted by $x_{n+1}$ and $\pd{t}$
by $\pd{n+1}$.

Let $L: \RR^{2n+2} \rightarrow \RR$ be a linear form
$L(\alpha,\beta)=\sum_i u_i \alpha_i + v_i \beta_i$ such that
$u_i+v_i\geq 0$ for $i=1,\ldots,n+1$, inducing the so-called
$L$--filtration on the ring $D$. If $u_i+v_i>0$ for all $i$,  the
associated graded ring $\gr^{L}(D)$ is isomorphic to a polynomial
ring in $2n+2$ variables
$(x,\xi)=(x,\xi_1,\ldots,\xi_{n+1})$ with complex
coefficients. This polynomial ring is $L$-graded, the $L$-degree of
a monomial $x^\alpha \xi^\beta$ being $L(\alpha,\beta)$. If we need
to emphasize the coefficients of the linear form we simply
write $L=L_{(u,v)}$ for $(u,v)\in \RR^{2n+2}$ with $u_i+v_i\geq 0$
for all $i$. If $u={\bf 0}\in \NN^{n+1}$ and $v={\bf
1}=(1,1,\ldots,1)\in \NN^{n+1}$ then the corresponding $L_{(u,v)}$
filtration is nothing but the usual order filtration on $D$ (which
is also called the $F$-filtration). If $u=({\bf 0},-1)\in \NN^{n+1}$
and $v=-u\in \NN^{n+1}$ then the corresponding $L_{(u,v)}$
filtration is nothing but the Malgrange-Kashiwara filtration on $D$
(also known as the $V$-filtration) with respect to $t=0$. In the
remainder of this section we assume  $u_i+v_i>0$ for all
$i$; we say then that $(u,v)$ is a {\em weight vector} for the Weyl
algebra $D$.

All the $D$--modules appearing here are left $D$-modules unless stated otherwise. We denote by
$T\subset \CC^{n+1}$ the hyperplane defined by $t=0$.

Let $M$ be a finitely generated $D$-module. To the $L$-filtration on
$D$ we associate a {\em good} $L$-filtration on $M$, by means of a
finite presentation. The associated $\gr^{L}(D)$-module $\gr^{L}(M)$
is then finitely generated. The radical of the annihilating ideal
$Ann_{\gr^{L}(D)}(\gr^{L}(M))$, which is independent of the {\em
good} $L$--filtration on $M$, defines an affine algebraic subset of
the cotangent space $T^*\CC^{n+1}=\CC^{2n+2}$.  This algebraic set
is called the $L$-{\em characteristic variety} of $M$ and it is denoted by
$Ch^{L}(M).$  The results stated so far are well
known in $D$-module theory and generalize
\cite{Bernstein70} which treats the case of the $F$--filtration in $D$. The case of a general  $L$--filtration has been studied for example in
\cite{Laumon} and, with more details,  in \cite[Section 3.2]{Laurent} in the micro-differential setting which is slightly
different from the one in this paper.
See also \cite[Section 2]{ACG-how-1996} for an equivalent treatment better
adapted to effective computations for modules on the Weyl algebra.

We consider a special type of $L$-filtration: For any
real number $r\in \RR_{\geq 0}$ we denote by $L_r$ either the linear
form $L_r=F+rV$ or the filtration on $D$ given by the
$(2n+2)$--dimensional weight vector $(0,\ldots,0,-r,1,\ldots,1,1+r)$ where
$-r$ is placed in the $(n+1)^{th}$-component. Here $F$ (resp. $V$)
stands for the order filtration on $D$ (resp. the
Malgrange-Kashiwara filtration with respect to $T$).

\begin{definition}\label{def:slope} {\rm \cite[Section 3.4]{Laurent}}
Let $M$ be a finitely generated $D$--module.
Consider the projection $\Pi : T^{\ast} \CC^{n+1} \longrightarrow T$ defined by $$\Pi (x_1 ,\ldots ,x_n ,t ,\xi_1
,\ldots ,\xi_n , \xi_t )=(x_1 ,\ldots ,x_n ,0).$$ For any real number
$r>0$, let $I_T^r (M)$ be the closure of the projection by $\Pi$ of
the irreducible components of the $L_r$--characteristic variety
$Ch^{L_r}M \subset T^*\CC^{n+1}$ that are not
$(F,V)$--bihomogeneous. The real number $s=r+1>1$ is said to be a
{\em slope} of $M$ along $T$ at $p\in T$ if and only if $p\in I_T^r (M)$.
\end{definition}

As proved by Y. Laurent, see \cite[Section 3.4]{Laurent}, any slope is a rational number and the set of slopes of $M$ is finite. Moreover, Y. Laurent also proved {\em loc. cit.}  that  $s=r+1$ is a slope of $M$ along $T$ at $p\in T$ if and only if, in a neighborhood of $\Pi^{-1}(p)$, $Ch^{L_{r'}}(M)$ is not locally constant for $r' \in (r -\epsilon , r +\epsilon )$
with $\epsilon >0$ small enough. The irregularity of a
holonomic system with respect to a smooth hypersurface (\cite[D\'ef.
6.3.1]{Mebkhout-positivite}) is deeply
related with the slopes of the system defined with
respect to the given hypersurface \cite{Laurent-Mebkhout-pa-pa}.

\section{On the irregularity of $A$--hypergeometric
systems}\label{section-irregularity-GKZ}

Recall that $A$ is a $d\times n$  integer matrix of rank $d$ whose
columns $a_1, \ldots,a_n$ generate $\ZZ^d$ as $\ZZ$--module. We
denote by $H_A(\beta)$ the hypergeometric ideal associated with $A$
and the parameter vector $\beta \in \CC^d$ \cite{gkz1989} and by
$M_A(\beta)$ the corresponding hypergeometric system (also known as
$GKZ$--system). This system is the quotient of $D_n:=
\CC[x_1,\ldots,x_n]\langle \partial_1,\ldots, \partial_n\rangle$,
the Weyl algebra of order $n$, by the left ideal $H_A(\beta)$. In
this section we denote $X=\CC^n$.

\subsection{Slopes of $A$--hypergeometric systems}
We recall in this subsection some results from \cite{SW}, where $A$ is
assumed to be pointed, i.e. the columns of $A$ lie in a open half-space defined by a
hyperplane passing through the origin in $\RR^d$.
These results will not be applied to our matrix $A$ but only to the matrix $\wA$
(see Subsection \ref{slopes-mhgs}), which we assume to be pointed throughout this article.

We denote by $a_i$ the $i$-th column of $A$ for $i=1,\ldots,n$. The
$L$-characteristic variety of the hypergeometric system $M_A
(\beta)$ has been described, in a combinatorial way, by M. Schulze
and U. Walther \cite{SW} for any pointed matrix $A$ and any
filtration $L=(u,v)$ such that $u_i + v_i =c>0$ for all $i=1,\ldots
,n$. The case $L=F$ was first studied  by A. Adolphson
\cite{Adolphson}. The main tool for their description is the notion of
$(A,L)$-{\em umbrella} that we define here for the sake of completeness. First of all, the $(A,L)$--umbrella only depends on $A$
and on the coefficients $v_i$ of the linear form $L=L_{(u,v)}$.

\begin{definition} {\rm \cite[Def. 2.7]{SW}} \label{special-umbrella}
We assume that $v_i>0$ for all $i$. The $(A,L)$--polyhedron $\Delta_A^{L}$ is the convex hull in $\RR^d$
of the set $\{{\bf{0}}, a_1/v_1,\ldots,a_n/v_n\}$. The
$(A,L)$-umbrella $\Phi_A^L$ is the set of faces of $\Delta_A^L$ which do not contain zero. In
particular, $\Phi_A^L$ contains the empty face.
\end{definition}

By $\Phi_A^{L,k} \subset \Phi_A^L$ we denote the subset of faces of
dimension $k$. We identify each face $\sigma$ of $\Phi_A^L$
with the set $\{ i: \; a_i /v_i \in \sigma \}$ and with $\{ a_i : \;
a_i /v_i \in \sigma \}$. The $(A,L)$-umbrella is then an abstract
cell complex.

When not all the $v_i$ are strictly positive then both
definitions of the $(A,L)$-polyhedron and the $(A,L)$-umbrella are a
little bit more involved. We refer to {\cite[Def. 2.7]{SW}} for
these precise definitions in the general case. See also Subsection
\ref{slopes-mhgs}.

\begin{theorem} {\rm \cite[Cor. 4.17]{SW}} \label{SW} The $L$--characteristic variety of $M_A (\beta)$ is given by:
\begin{equation}
 Ch^L (M_A (\beta)) =\bigcup_{\tau \in \Phi_A^L} \overline{C_A^{\tau}}
\end{equation} where $\overline{C_A^{\tau}}$ is the closure of the conormal $C_A^{\tau}$ of the torus orbit
$O_A^{\tau}=\{\xi \in T_{0}^{\ast} X =\CC^n : \;
\xi_i = 0 \mbox{ if } i\notin \tau , \; \xi_i = y^{a_i} \; \mbox{ if } i\in \tau , y \in (\CC^{\ast})^d   \}$.
\end{theorem}

Moreover, it is proved in \cite[Lemma 3.14]{SW} that
$\overline{C_A^{\tau}}$ meets $T_{0}^{\ast}X$  for all $\tau \in
\Phi_A^{L}$. Thus Theorem \ref{SW} provides a description of the
slopes of a hypergeometric system at the origin along any coordinate
variety $Y\subset X$ via considering the $1$--parameter
family of filtrations $L_r = F + r V$, $r>0$, where $V$ is the $V$-filtration\footnote{The $V$--filtration with respect to the coordinate variety $Y=(x_1=\cdots=x_\ell=0)$ is defined
by assigning the weight -1 (resp. the weight 1) to the variables $x_i$ (resp. $\partial_i$) for $i=1,\ldots,\ell$ and the weight $0$ to the remaining variables.} along $Y$:

\begin{corollary} {\rm \cite[Cor. 4.18]{SW}} \label{SW2}
The real number $s=r+1>1$  is a slope of $M_A (\beta )$ along $Y$ at
the origin if and only if $\Phi_A^{L_{r'}}$ is not locally constant at $r'=r$.
\end{corollary}

\begin{remark}\label{remark-slope-at-any-p}
When $Y\subset X$ is a coordinate hyperplane then the set of slopes of $M_A
(\beta)$ along $Y$ at the origin coincides with the set of slopes of $M_A (\beta)$ along $Y$ at any point $p\in Y$. This is proved in \cite[Th. 5.9]{Fer} by using the
comparison result in \cite[Th. 2.4.2]{Laurent-Mebkhout-pa-pa}.
\end{remark}

\subsection{Gevrey solutions of $A$--hypergeometric systems at infinity.}\label{Gevrey-solutions-at-infinity}

In this subsection we extend some of the results from \cite{Fer}
to the case of the Gevrey solutions of a hypergeometric system at
infinity in the direction of a coordinate hyperplane that we may assume to be $x_n=0$.  This
construction is used later in the study of the Gevrey solutions
along $T$ of a modified hypergeometric system, see Subsection \ref{Gevrey-solutions-modulo-convergent-series}.

Let us denote by $\cO_{X}$ the sheaf of holomorphic functions on
$X=\CC^n$. For $Y=\{x_n = 0\}$, we denote by $\cO_{\widehat{X\vert
Y}}$ the formal completion of $\cO_{X}$ along $Y$, whose germs at
$(p,0)\in Y$ are of the form $f=\sum_{m\geq 0} f_m x_n^m$ where  all
the $f_m=f_m (x_1 ,\ldots , x_{n-1})$ are holomorphic functions in a
common neighborhood of $p$. Notice that the restriction of $\cO_X$
to $Y$, denoted by $\cO_{X\vert Y}$, is a subsheaf of
$\cO_{\widehat{X\vert Y}}$. For any real number $s$, we also
consider the sheaf $\cO_{\widehat{X\vert Y}}(s)$ of Gevrey series
along $Y$ of order $s$ which is defined to be the subsheaf of $\cO_{\widehat{X\vert Y}}$
whose germs $f$ at $(p,0)\in Y$ satisfy $$\sum_{m\geq 0}
\frac{f_m}{{m!}^{s-1}} x_n^m  \in \cO_{X \vert Y, (p,0)}.$$
We denote $\cQ_Y(s):=\frac{\cO_{\widehat{X\vert Y}}(s)}{\cO_{X\vert
Y}}$ and use $\cO_{\widehat{X\vert Y}}({<}\,\!s)$ for the sheaf of
Gevrey series along $Y$ of order less than $s$. If $f$ belongs to
$\cO_{\widehat{X\vert Y}}(s) \setminus \cO_{\widehat{X\vert
Y}}({<}\,\!s)$ for some $s$, we say that the {\em index} of the
Gevrey series $f$ is $s$. We also write $\cO_{\widehat{X\vert
Y}}(+\infty) := \cO_{\widehat{X\vert Y}}$.

We  denote by $\cD=\cD_X$ the sheaf of linear differential
operators on $X$ with holomorphic coefficients. If $\cM$ is a
coherent $\cD$-module, Z. Mebkhout has defined in \cite[D\'ef.
6.3.1]{Mebkhout-positivite} the irregularity of order $s$ of $\cM$
with respect to $Y$ to be  the complex of (sheaves of) vector spaces
$Irr^{(s)}_Y(\cM):=\RR \cH om_\cD(\cM, \cQ_Y(s))$ and has proved
that, for all $s\in [1,+\infty]$,  this complex is a perverse sheaf
on $Y$  when  $\cM$ is holonomic \cite[Th.
6.3.3]{Mebkhout-positivite}. By the comparison theorem
\cite[Th.2.4.2]{Laurent-Mebkhout-pa-pa} $s$ is a slope of $\cM$ with
respect to $Y$ if and only if $s$ is a gap in the filtration
$Irr^{(s)}_Y(\cM)$ on the irregularity
$Irr_Y(\cM):=Irr_Y^{(+\infty)}(\cM)$. The perversity result implies,
in particular, that  at a generic point $p\in Y$ only the first
cohomology space of previous complexes is possibly non zero and thus it is worth studying the stalk $\cH om_\cD(\cM, \cQ_Y(s))_p$.

Recall that $A=(a_1 \; \cdots \; a_n )$ is a full rank matrix with
$a_i \in \ZZ^d$ for all $i=1 ,\ldots ,n$ and $d\leq n$. Following
\cite{gkz1989} and \cite{SST}, for any vector $v\in \CC^n$
we can define a series
\begin{equation}
\phi_v=\phi_v(x):=\sum_{u\in N_v} \frac{[v]_{u_{-}}}{[v+u]_{u_{+}}}
x^{v+u} \label{phiv}\end{equation} where $v\in \CC^n$ verifies $A
v=\beta$ and $N_v = \{u\in \ker (A)\cap \ZZ^n  : \;
\operatorname{nsupp}(v+u)=\operatorname{nsupp}(v) \}$. Here $\ker
(A) =\{u\in \mathbb{Q}^n : \; A u=0 \}$,
$\operatorname{nsupp}(w):=\{i\in \{1,\ldots , n\}: \; w_i \in
\ZZ_{<0 }\}$ is the negative support of $w\in \CC^n$,
$[v]_{u}=\prod_{i} [v_i ]_{u_i} $ and $[v_i ]_{u_i}=\prod_{j=1}^{u_i
} (v_i -j +1)$ is the Pochhammer symbol for $v_i \in \CC$, $u_i \in
\NN$.

The series $\phi_v$ is annihilated by the hypergeometric ideal $H_A
(\beta )$ if and only if the negative support of $v$ is minimal,
i.e., $\nexists u \in \ker (A)\cap \ZZ^n $ with
$\operatorname{nsupp}(v+u)\subsetneq \operatorname{nsupp}(v)$ (see
\cite[Section 3.4]{SST}).

When $\beta \in \CC^d$ is {\em very generic}, i.e., when $\beta $ is not
in a locally finite countable union of Zarisky closed sets, there is a basis of the
Gevrey solution space of $\cM_A (\beta )$ along $Y$ given by series
$\phi_v$ for suitable vectors $v\in \CC^n$ (see \cite[Ths. 6.2 and
6.7]{Fer}).

For any subset $\eta\subseteq \{1,\ldots ,n\}$ we denote by
$A_{\eta}$ the submatrix of $A$ given by the columns of $A$ indexed
by $\eta$ and we denote $\bar \eta =\{1,\ldots ,n\}\setminus \eta$.

We say that $\sigma \subseteq \{1,\ldots ,n \}$ is a $(d-1)$-{\em simplex}
with respect to $A$ (or simply that $\sigma$ is a $(d-1)$-simplex) if
the columns of $A_\sigma$ determine a basis of $\RR^d$. If it is so,
we can reorder the variables in order to have $\sigma=\{1,\dots
,d\}$ without loss of generality. Then a basis of $\ker (A)$
associated with $\sigma$ is given by the columns of the matrix:
$$B_{\sigma}=
\left(\begin{array}{cccc}
               -A_{\sigma}^{-1} a_{d+1} & -A_{\sigma}^{-1}a_{d+2} & \cdots & -A_{\sigma}^{-1}a_{n} \\
               1 & 0 &  & 0 \\
               0 & 1 &  & 0 \\
          \vdots &   & \ddots & \vdots \\
               0 & 0 &   & 1
             \end{array}\right)$$

A vector $v\in \CC^n$ satisfying $v_i \in \NN$ for all $i\notin
\sigma$ and $A v=\beta$ can be written as
$$v^{\mathbf{k}}=(A_{\sigma}^{-1}(\beta - \sum_{i\notin \sigma} k_i
a_i),\mathbf{k})$$ for some  $ \mathbf{k}=(k_i )_{i\notin \sigma
}\in \NN^{n-d} $. Since $\beta$ is very generic then the negative
support of $v^{\mathbf{k}}$ is the empty set and hence
$\phi_{v^{\mathbf{k}}}$ is annihilated by $H_A (\beta )$. Moreover,
the summation index set $N_{v^{\mathbf{k}}}$ in the series
$\phi_{v^{\mathbf{k}}}$ is given by the integer vectors in an affine
translate of the positive span of the columns of $B_{\sigma}$. The
sum of the coordinates of the $i$-th column of $B_{\sigma}$ is $
1-|A_{\sigma}^{-1}a_{d+i}|$ where $| \mbox{ }|$ means the sum of the coordinates.
We have the following.

\begin{theorem}{\rm \cite[Theorem 3.11]{Fer}}\label{Theorem-Gevrey-solutions-hypergeometric}
Under the above conditions the series $\phi_{v^{\mathbf{k}}}$ is a
Gevrey solution of $M_A (\beta)$ along
$Z=\{x_j =0 :\; |A_{\sigma}^{-1}a_{j}| >1\}$
with index $s=\max \{|A_{\sigma}^{-1}a_{j}|: j=1 ,\ldots ,n\}$ at points in certain
relatively open subset of $Z$. In particular, if $|A_{\sigma}^{-1}a_{j}|\leq 1$
for all $1\leq j\leq n$ then $\phi_{v^{\mathbf{k}}}$ is convergent.
\end{theorem}

By Corollary \ref{SW2}, a real number  $s>1$ is a slope of $M_A(\beta
)$ along $Y=\{x_{n} =0\}$ if and only if  $\frac{1}{s}a_{n}$ belongs
to the hyperplane $H_{\tau}$ supported on a facet $\tau$ of the
convex hull of $\{0, a_{1},\ldots ,a_{n-1}\}$ such that $0\notin
\tau$. In particular, for any $(d-1)$-simplex $\sigma \subseteq \tau$ it
is easy to check that $s=|A_{\sigma}^{-1}a_n|>1$. We say in this
case that $\sigma$ is a $(d-1)$-simplex corresponding to the slope $s>1$
of $M_A (\beta )$ along $Y$.

The following Theorem is a summary of some of the results from \cite{Fer}.
Its last statement uses results from \cite{Laurent-Mebkhout-pa-pa} and in \cite{SW}. For the definition
of a regular triangulation see, e.g., \cite[Ch. 8]{Sturmfels}.

\begin{theorem}\label{hypergeometric-Gevrey-solutions}
Assume that $\beta \in \CC^d$ is very generic and that $s>1$ is a
slope of $M_A (\beta )$ along $Y=\{x_{n} =0\}$. For any $(d-1)$-simplex
$\sigma$ corresponding to $s$ one can construct
$\operatorname{vol}(\sigma)=|\det (A_{\sigma})|$ many linearly
independent Gevrey solutions $\phi_{v^{\mathbf{k}}}$ of $M_A (\beta
)$ along $Y$ with index $s$ by varying $ \mathbf{k} \in \NN^{n-d}$
in a set $\Lambda$ so that $\{A_{\overline{\sigma}}\mathbf{k}\, \vert
\mathbf{k} \in  \Lambda \}$ is a set of representatives of the group
$\ZZ^d /\ZZ A_{\sigma}$. Moreover, if we repeat this construction
for all the $(d-1)$-simplices $\sigma $ corresponding to $s$ which
belong to a suitable regular triangulation for $A$ and take the
classes modulo $\cO_{\widehat{X|Y}}({<}\,\!s)$ we obtain a basis for
the space of solutions of $M_A (\beta )$ in the space
$(\cO_{\widehat{X|Y}} (s)/\cO_{\widehat{X|Y}}({<}\,\!s))_p$ for
points $p\in Y$ in a relatively open set of $Y$.
\end{theorem}

We have exhibited the construction of the Gevrey solutions of $M_A
(\beta )$ along $Y=\{ x_n = 0\}$ corresponding to each slope $s>1$
of $M_A (\beta )$ along $Y$ for $\beta$ very generic.

Let us construct Gevrey solutions of $M_A (\beta )$ at infinity. In
other words, we are going to construct Gevrey solutions of the
projectivized hypergeometric system treated in \cite[Section 5]{SW2}
at a generic point at infinity. We use the following notation:
$X'$ is $\CC^{n}$ with coordinates $(x_1 ,\ldots ,x_{n-1}, z)$
and $z=1/x_n$ so that $X\cap X'=\CC^{n-1} \times \CC^{\ast}$. Denote
$Y'=\{x_n=\infty\}=\{z=0\}\subseteq X'$.

Take $L_{-r}=F-r V$ where $V$ is the $V$-filtration along $Y$.
Notice that $\Phi_A^{L_{-r}}$ is not locally constant at $r=s-1>0$
if and only if $\frac{1}{(2-s)}a_n$ belongs to the hyperplane
$H_{\tau}$ supported on a facet $\tau$ of the convex
hull of $\{0, a_{1},\ldots ,a_{n-1}\}$ such that $0\notin \tau$.

\begin{theorem}
\label{hypergeometric-Gevrey-solutions-infinity} Assume that $\beta
\in \CC^d$ is very generic and that there exists $s>1$ such that
$\frac{1}{(2-s)}a_n$ belongs to a hyperplane $H_{\tau}$ as above.
For any $(d-1)$-simplex $\sigma\subseteq \tau $ one can construct
$\operatorname{vol}(\sigma)=|\det (A_{\sigma})|$ many linearly
independent Gevrey series $\phi_{v^{\mathbf{k}}}$ along $Y'$ with
index $s$ by varying $ \mathbf{k} \in \NN^{n-d-1} \times \ZZ_{<0}$
in a set $\Lambda$ so that $\{A_{\overline{\sigma}}\mathbf{k} :\; \mathbf{k} \in  \Lambda \}$ is a
set of representatives of the group
$\ZZ^d / \ZZ A_{\sigma}$.
The classes of these series modulo convergent series $\cO_{X'|Y'}$ are
solutions of $M_A (\beta)$ in $\cO_{\widehat{X'|Y'}} (s) /\cO_{X'|Y'}$.

Moreover, if we repeat this construction for all the $(d-1)$-simplices
$\sigma $ as above which belong to a suitable regular triangulation
for $A$ and take the classes modulo $\cO_{\widehat{X'|Y'}}({<}\,\!s)$ then
we obtain a basis for the space of solutions of $M_A (\beta )$ in
the space $(\cO_{\widehat{X'|Y'}} (s)/\cO_{\widehat{X'|Y'}}
({<}\,\!s))_p$ for points $p\in Y'$ in a relatively open set of $Y'$.
\end{theorem}

{\it Proof.-}
By reordering the variables we may  assume
$\sigma=\{1,\dots ,d\}$ without loss of generality. Let
$B_{\sigma}'$ be the matrix $B_{\sigma}$ but with the last column
multiplied by $-1$.

Take a vector $v\in \CC^n$ such that $A v=\beta$, $v_i \in \NN$ for
all $i \notin \sigma \cup \{n\}$ and $v_n \in \ZZ_{<0}$. It is clear
that such a vector can be taken as $v=v^{\mathbf{k}}$, as before,
for some  $ \mathbf{k}=(k_i )_{i\notin \sigma }\in \NN^{n-d-1}\times
\ZZ_{<0}$. Moreover, the summation index set $N_{v^{\mathbf{k}}}$ in
the series $\phi_{v^{\mathbf{k}}}$ is given by the integer vectors
in an affine translate of the positive span of the columns of
$B_{\sigma} '$. Notice that the sum of the coordinates of the $i$-th
column of $B_{\sigma} '$ is $ 1-|A_{\sigma}^{-1}a_{d+i}|\geq 0$ for
$i=1,\ldots , n-d -1$ while the sum of the coordinates of the last
column is $|A_{\sigma}^{-1}a_n |-1 =2-s -1 =1-s <0$. This implies
that
$\phi'_{v^{\mathbf{k}}}(x_1,\ldots,x_{n-1},z):=\phi_{v^{\mathbf{k}}}(x_1,\ldots,x_{n-1},1/z)$
is a Gevrey series along $Y'$ with order $s=r+1>1$. To prove that
$s$ is the Gevrey index of $\phi'_{v^{\mathbf{k}}}$ one can use a
slightly modified version of Lemma 3.8 in \cite{Fer}.

Moreover, since $\beta$ is very generic then the negative support of
$v^{\mathbf{k}}$ is $\operatorname{nsupp}{(v^{\mathbf{k}})}=\{n\}$
which is not minimal and hence $\phi_{v^{\mathbf{k}}}$ is not
annihilated by $H_A (\beta )$ (see \cite[Section 3.4]{SST}).
However, it can be checked that for any differential operator $P\in
H_A (\beta )$ the series $P(\phi_{v^{\mathbf{k}}})$ is either zero
or a polynomial in $z = x_n^{-1}$ with coefficients that are
convergent power series in the variables $x_1 ,\ldots , x_{n-1}$.
Now we finish the proof following \cite[Sec. 7]{Fer} and  \cite[Sec. 5]{SW2}.  \QED

\begin{remark}\label{remark-Gevrey-solutions-infinity}
A slightly weaker version of the first paragraph of Theorem
\ref{hypergeometric-Gevrey-solutions-infinity} can also be proved
without the very genericity assumption in $\beta \in \CC^d$. More
precisely, for all $\beta \in \CC^d$ if $\Phi_A^{L_{-r}}$ is not
locally constant at $r=s-1>0$ we can construct a Gevrey series along
$Y'$ of index $s=r+1$ which is a solution of $M_A (\beta )$ modulo
the space $\cO_{\widehat{X'|Y'}}({<}\,\!s)$ (of Gevrey series along
$Y'$ with order less than $s$) by methods similar to the ones in
\cite[Section 4]{Fer}.
\end{remark}

The following corollary is a particular case of \cite[Conjecture 5.18]{SW2}.

\begin{corollary}\label{slopes-hypergeometric-at-infinity}
The real number $s>1$ is a slope of $M_A (\beta )$ along $Y'=\{x_n =\infty \}$
if and only if $\Phi_A^{L_{-r}}$ is not locally constant at $r=s-1$.
\end{corollary}

{\it Proof.-} For the only if direction of the proof
we refer to \cite[Section 5]{SW2}. Let us prove the if direction. By Theorem
\ref{hypergeometric-Gevrey-solutions-infinity} and Remark
\ref{remark-Gevrey-solutions-infinity} one can construct a Gevrey
series of index $s=r+1$ that is a solution of $M_A (\beta)$ along $Y'$ (modulo $\cO_{\widehat{X'\vert Y'}}({<}\,\!s)$).
So, the result follows from the comparison theorem for the slopes \cite{Laurent-Mebkhout-pa-pa}.
\QED

\begin{example}\label{example-slope-infinity}
Set $$A =
\left( \begin{array}{ccc}
            2 & 1 & 1 \\
            1 & 2 & 1 \end{array} \right)$$ and $\beta \in \CC^2$. We have that the kernel of $A$ is
generated by $u=(1,1,-3)$ and thus
the hypergeometric ideal $H_A (\beta )$ is generated by
$\Box_u = \partial_1 \partial_2 - \partial_3^3$,
$E_1 -\beta_1 =  2 x_1 \partial_1 + x_2 \partial_2  + x_3 \partial_3 -\beta_1$ and
$E_2 -\beta_2 =  x_1 \partial_1 + 2 x_2 \partial_2  + x_3 \partial_3  -\beta_2 $.
Take $\sigma =\{1 ,2 \}$ and notice that for $s=4/3$ we have that $a_3 / (2-s)$ belongs to the line $H_{\sigma}$
determined by $a_1$ and $a_2$.
By Corollary \ref{slopes-hypergeometric-at-infinity} we have that $s=4/3$ is a slope of $M_A (\beta )$ along
$Y'=\{x_3 = \infty\}$. Indeed, if we consider $v^{\mathbf k}=
((2 \beta_1 -\beta_2 -k )/3 ,(2 \beta_2 -\beta_1 - k)/3 , k)$ for
${\mathbf k}=k\in \ZZ_{< 0}$ we have that $\phi_{v^{\mathbf k}}$ is a Gevrey series along $Y'$ of index
$s=4/3$ if $v_1^{\mathbf k} , v_2^{\mathbf k} \notin \ZZ_{<0}$. Moreover,
we have  that $(E_i -\beta_i)( \phi_{v^{\mathbf k}})=0$ for $i=1,2$.
For each of the three biggest $k\in \ZZ_{< 0}$
verifying that $v_1^k , v_2^k \notin \ZZ_{<0}$ we have that
$\Box_u (\phi_{v^k})= v_1^k v_2^k x_1^{v_1^k -1} x_2^{v_2^k -1} x_3^{k}$,
which is convergent at any $p\in Y' \cap \{ x_1 x_2 \neq 0\}$. The classes modulo
$\cO_{X'|Y'}$ of these three series $\phi_{v^{\mathbf k}}$ form a basis for the space of solutions of
$M_A (\beta )$ in $(\cO_{\widehat{X'|Y'}}(s)/\cO_{X'|Y'})_p$, $p\in Y' \cap \{ x_1 x_2 \neq 0\}$.
If $\beta$ is very generic then $k=-1,-2,-3$.
\end{example}

\section{On the irregularity of modified $A$-hypergeometric
systems}\label{section:irregularity-MGKZ}

\subsection{Fourier transform and initial
ideals}\label{section-Fourier} Let $D$ be the Weyl algebra
${\CC}\langle x_1, \ldots, x_n, t,\pd{1}, \ldots, \pd{n}, \pd{t}
\rangle$. The variable $t$ is also denoted by $x_{n+1}$ and $\pd{t}$
by $\pd{n+1}$.

Let $L: \RR^{2n+2} \rightarrow \RR$ be a linear form
$L(\alpha,\beta)=\sum_i u_i \alpha_i + v_i \beta_i$ such that
$u_i+v_i\geq 0$ for $i=1,\ldots,n+1$ inducing the so-called
$L$--filtration on the ring $D$. If $u_i+v_i>0$ for all $i$,  the
associated graded ring $\gr^{L}(D)$ is isomorphic to a polynomial
ring in $2n+2$ variables
$(x,\xi)=(x_1,\ldots,x_{n+1},\xi_1,\ldots,\xi_{n+1})$ with complex
coefficients. This polynomial ring is $L$-graded, the $L$-degree of
a monomial $x^\alpha \xi^\beta$ being $L(\alpha,\beta)$. If we need
to emphasize the coefficients of the linear form we will simply
write $L=L_{(u,v)}$ for $(u,v)\in \RR^{2n+2}$ with $u_i+v_i\geq 0$
for all $i$. If $u={\bf 0}\in \NN^{n+1}$ and $v={\bf
1}=(1,1,\ldots,1)\in \NN^{n+1}$ then the corresponding $L_{(u,v)}$
filtration is nothing but the usual order filtration on $D$ (which
is also called the $F$-filtration). If $u=({\bf 0},-1)\in \NN^{n+1}$
and $v=-u\in \NN^{n+1}$ then the corresponding $L_{(u,v)}$
filtration is nothing but the Malgrange-Kashiwara filtration on $D$
(we also say the $V$-filtration) with respect to $t=0$. In the
remainder of this section we will assume  $u_i+v_i>0$ for all
$i$; we will say then that $(u,v)$ is a weight vector for the Weyl
algebra.

We define the ring isomorphism $\F$ of $D$ by $ t
\mapsto -\pd{t}$, $ \pd{t} \mapsto t$. The isomorphism $\F$ is
called the {\it Fourier transform} on $D$ with respect to the
variable $t$. The inverse transform $\F^{-1}$ is given by $ t
\mapsto \pd{t}$, $ \pd{t} \mapsto -t$. Let $(u,v)$ be a weight
vector for the Weyl algebra. The Fourier transforms $\F$ and
$\F^{-1}$ induce isomorphisms in ${\rm gr}^L(D)$ and we denote them
also by $\F$ and $\F^{-1}$ respectively. Analogously, if $C \subset
\CC^{2n+2}$ is the affine algebraic set defined by an ideal
$J\subseteq {\rm gr}^L(D)$ we write $\F C$ and $\F^{-1}C$ for the
algebraic set defined by $\F J \subseteq {\rm gr}^L(D)$ and $\F^{-1}
J \subseteq {\rm gr}^L(D)$ respectively.

We define the Fourier transform of the weight vector $(u,v)$ by the
formula $$\F(u,v):=(u_1, \ldots, u_n, v_{n+1}, v_1, \ldots, v_n,
u_{n+1}).$$ We notice that $\F \F (u,v)=(u,v)$. We will also write
$\F L=\F(u,v)$ if $L=L_{(u,v)}$.

\begin{proposition}\label{FinF(u,v)FP}
For any operator $\ell \in D$, we have
$${\rm in}_{(u,v)}(\ell) = \F^{-1} {\rm in}_{\F (u,v)} (\F \ell) $$
\end{proposition}

{\it Proof}\/. We prove it in the case $n=0$. Other cases can be
reduced to this case.
We put $\xi_{n+1} = {\rm in}_{(u,v)}(\pd{t})$. We assume that $u+v
>0$ and $\ell = t^a \pd{t}^b$. Then, we have ${\rm in}_{(u,v)}(\ell)
= t^a \xi_{n+1}^b$. Since $\F \ell = (-\pd{t})^a t^b = (-1)^a(t^b
\pd{t}^a + ab t^{b-1} \pd{t}^{a-1} + \cdots)$ and $u+v > 0$, we have
${\rm in}_{(v,u)}(\F \ell) = (-1)^a t^b \xi_{n+1}^a$. Applying the
inverse Fourier transform, we obtain $\F^{-1} {\rm in}_{(v,u)}(\F
\ell) = (-1)^a \xi_{n+1}^b (-t)^a = {\rm in}_{(u,v)}(\ell)$.

Suppose that $t^a \pd{t}^b >_{(u,v)} t^{a'} \pd{t}^{b'}$. Then, we
have $(-\pd{t})^a t^b >_{(v,u)} (-\pd{t})^{a'} t^{b'}$. Thus, we
obtain the conclusion. \QED

Proposition \ref{FinF(u,v)FP} yields the following simple, but
important claim for the Gr\"obner deformation method.

\begin{corollary} \label{corollary:fourier}
For any left ideal $I$ in  $D$, we have
$$ {\rm in}_{(u,v)}(I) = \F^{-1} {\rm in}_{\F(u,v)}(\F I)$$
\end{corollary}

{\it Proof}\/. Since ${\rm in}_{(u,v)}(I)$ is spanned by ${\rm
in}_{(u,v)}(\ell)$, $\ell \in I$ as a ${\CC}$-vector space, the
conclusion follows from the previous proposition. \QED

\subsection{Slopes of modified $A$-hypergeometric systems}\label{slopes-mhgs}

We retain the notations of \cite{takayama2009}, which are explained
in the introduction. We are interested in the slopes of the modified
system $M_{A,w ,\alpha}(\beta)$ at any point along $T=\{t=0\}$.
We notice that $H_{A,w,\alpha}(\beta)$ is the Fourier transform of
$H_{\wA}(\tilde \beta)$ where ${\tilde \beta} = (\beta,\alpha -1)$,
i.e.
\begin{equation}\label{FouriermGKZ=GKZ}  M_{A,w,\alpha}(\beta) = \cF M_{{\widetilde A}(w)}(\widetilde{\beta}).
\end{equation}
Recall that $L_r = F+ r V=(0,\ldots,0,-r,1,\ldots,1,1+r)$ where
$V=(0, \ldots, 0,-1; 0, \ldots, 0, 1)$ and $F=({\bf 0},{\bf 1})$.
Thus $\cF L_r = (0,\ldots,0,1+r,1,\ldots,1,-r)$.

In order to obtain the slopes of the modified system
$M_{A,w,\alpha}(\beta)$, we need to study the $L_r$--initial ideal
of the modified hypergeometric ideal $H_{A,w,\alpha}(\beta)$, for
$r\in {\mathbb R}_{>0}$.  Therefore, applying Corollary
\ref{corollary:fourier} and (\ref{FouriermGKZ=GKZ}), we have
\begin{equation}
 \label{fourier2}
 {\rm in}_{L_r}(H_{A,w,\alpha}(\beta))
 = \F^{-1} {\rm in}_{\F L_r} (H_{\wA}(\tilde \beta)).
\end{equation}

Using (\ref{fourier2}), we have
\begin{equation}\label{characteristic-fourier}
 Ch^{L_r}(M_{A,w,\alpha}(\beta)) = \cF^{-1} (Ch^{\cF L_r}
(M_{\wA}(\widetilde{\beta})))
\end{equation}

On the other hand, since $M_{\wA}(\widetilde{\beta})$ is a
hypergeometric system and the matrix $\wA$ is pointed, Theorem
\ref{SW} gives a description of the irreducible components of
$Ch^{\cF L_{r}}(M_{\wA}(\widetilde{\beta}))$ in terms
of the $(\wA , \cF L_r )$--umbrella. We notice here that the last
coordinate of $\cF L_r$ equals $-r<0$. We recall the definition
\cite[Def. 2.7]{SW} of the umbrella in this case. First of all let
us recall {\it loc. cit.,} that if $a$, $b$ are two points and $H$
is a hyperplane in $\PP^{d+1}(\RR)$ containing neither $a$ nor $b$,
then the convex hull ${\rm conv}_{H}(a,b)$ of $a$ and $b$ relative
to $H$, is the unique line segment joining $a$ and $b$ and not meeting $H$.

Let us denote the $i^{\rm{th}}$ column of $\wA$ by $\wa_i $, say $\wa_i
=\left(\begin{smallmatrix} a_i \\ w_i
\end{smallmatrix}\right)$ for $i=1,\ldots ,n$ and
$\wa_{n+1}=(0,\ldots,0,1)^t$. For simplicity let us write
$(\widetilde{A},\widetilde{L})$ instead of $(\wA , \cF L_r )$.  We
view $\wa_1,\ldots,\wa_{n+1}$ as points in ${\mathbb R}^{d+1}
\subset {\mathbb P}^{d+1}({\mathbb R})$. As $\widetilde{A}$ is
pointed, there exists a linear form $h$ on ${\mathbb R}^{d+1}$ such
that $h(\wa_i)>0$ for all $i$. Let $\epsilon \in {\mathbb R}$ be
such that $0<\epsilon <h(\wa_i)$ for $i=1,\ldots,n$ and $0<\epsilon
<\frac{h(\wa_{n+1})}{r}$.

\begin{definition} {\rm \cite[Def. 2.7]{SW}} \label{umbrella}
The $(\widetilde{A},\widetilde{L})$--polyhedron
$\Delta_{\widetilde{A}}^{\widetilde{L}}$ is the convex hull
$$\Delta_{\widetilde{A}}^{\widetilde{L}}= {\rm
conv_{H_\epsilon}}(\{{\bf{0}}, \wa_1,\ldots,\wa_n,
\frac{\wa_{n+1}}{-r}\})\subset \PP^{d+1}(\RR)$$ where $H_\epsilon$
is the projective closure of the affine hyperplane
$h^{-1}(-\epsilon)$. The $(\widetilde{A},\widetilde{L})$-umbrella
$\Phi_{\widetilde{A}}^{\widetilde{L}}$ is the set of faces of
$\Delta_{\widetilde{A}}^{\widetilde{L}}$ which do not contain the
origin. In particular, $\Phi_{\widetilde{A}}^{\widetilde{L}}$
contains the empty face. \end{definition}

{\begin{figure}[h!]\centering
\fbox{\includegraphics*[scale=1]{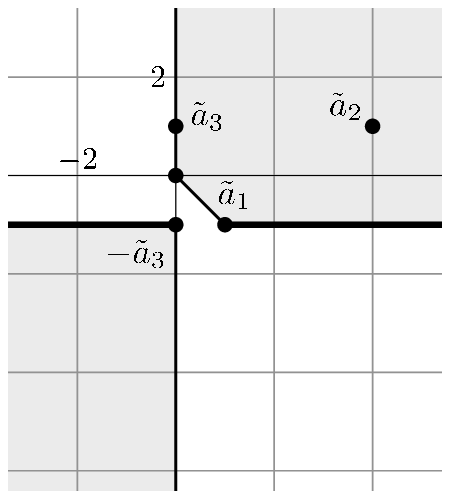}}
\fbox{\includegraphics*[scale=1]{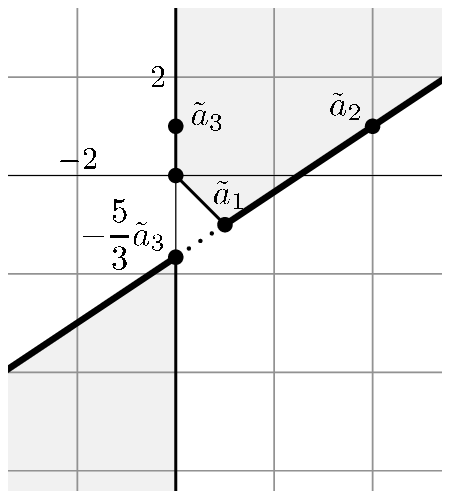}}
\fbox{\includegraphics*[scale=1]{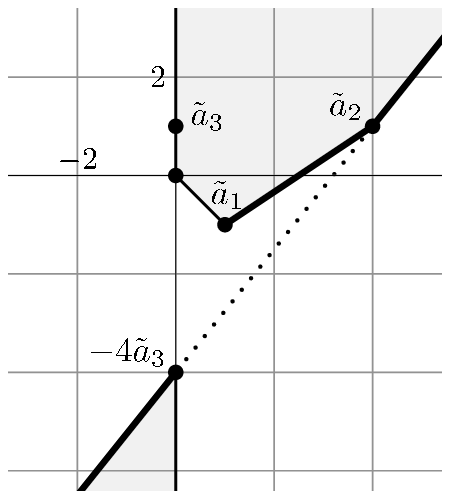}}
\caption{The first umbrella is for $r=1$, the second is for
$r=3/5$, and the third is for $r=1/4.$}
\end{figure}}

Figure 1 shows three  $(\wA , \cF L_r )$--umbrellas for $A=(1,4)$ and
$w=(-1,1)$. In each case the shaded region is the polyhedron
$\Delta_{\widetilde{A}}^{\widetilde{L}}$. For
$r=3/5$ the point $\frac{\widetilde{a}_3}{-r}$ belongs to the line
passing through $\widetilde{a}_1$ and $\widetilde{a}_2$ which  means
that $s=r+1=8/5$ is a slope of the system along $x_3=0$.

Using Definition \ref{def:slope}, equation
\eqref{characteristic-fourier} and \cite[Cor. 4.12]{SW} we get the
following:

\begin{corollary} \label{FourierSW} The real number $s=r+1>1$ is a slope of $M_{A,w,\alpha}(\beta)$ along $T$ at $p\in T$
if and only if $\Phi_{{\widetilde A}(w)}^{\cF L_{r'}}$ is not locally constant
at $r'=r$.
\end{corollary}

Let us denote by $A_w$ the matrix with columns $\wa_i$, $1\leq i\leq
n$, and let $\Delta_{A_w}$ be the convex hull of $\{\wa_i :\; 1\leq
i \leq n \}$ and the origin. With this notation, Corollary
\ref{FourierSW} can be rephrased as follows.

\begin{corollary} \label{FourierSW2} The real number $s=r+1>1$ is a slope of $M_{A,w, \alpha}(\beta)$
along $T$ at $p\in T$ if and only if
there exists a facet $\tau$ of $\Delta_{A_w }$ such that $0\notin
\tau$ and $-\dfrac{1}{r} \wa_{n+1}\in H_{\tau}$, where $H_{\tau}$ is
the hyperplane that contains $\tau$.
\end{corollary}

\begin{remark}\label{Remark-slopes-modify-infinity}
We note that the $D$-modules $M_{A,w,\alpha}(\beta)$ and
$M_{A,-w,-\alpha}(\beta)$ agree on $\CC^n \times \CC^{*}$ under the
change of the variable $t \rightarrow 1/t$ in $t\neq 0$. We regard
both $D$-modules as extensions of each other. In this paper, when we
say the irregularity of $M_{A,-w,-\alpha}(\beta)$ along
$T'=\{t=\infty\}$, it means the irregularity of
$M_{A,w,\alpha}(\beta)$ along $T=\{t=0\}$. Using this terminology,
Corollary \ref{FourierSW2} provides a description of the slopes of
the modified system $M_{A,w,\alpha}(\beta)$ along $T'$ by using
$\Delta_{A_{-w} }$ instead of $\Delta_{A_w }$.
\end{remark}

Until the end of this section we will denote either
$\Phi_{\widetilde{A}}^{v}$ or $\Phi_{\widetilde{A}}^{L}$ for the
$(\widetilde{A},L)$--umbrella with $L=L_{(u,v)}$ since Definition
\ref{umbrella} does not depend on $L$ but only on $v \in \RR^{n+1}$.
Moreover, for any subset $\eta\subseteq \{1,\ldots , n\}$ we denote
by $w_\eta$ the vector with coordinates equal to the ones of $w$
indexed by $\eta$, i.e., $w_\eta= (w_i )_{i\in \eta}$.

In the following two lemmas we assume for simplicity that $w_i >0$, $i=1,\ldots,n$.
In fact, this can be assumed without loss of generality since
$\wA $ is pointed.

\begin{lemma}\label{lemma-umbrellas-1}
For any sufficiently small real number $r>0$, we have that
$$\{\eta ' \in \Phi_{\wA}^{\cF L_r , d} :\; n+1 \in \eta ' \} = \{ \sigma \cup \{ n+1 \}: \; \sigma \in \Phi_{A_{\eta}}^{w_{\eta}, d-1},  \eta \in
\Phi_A^{F, d-1}\}.$$ In particular, if all the facets in
$\Phi_A^{F}$ contain exactly $d$ columns of $A$, then the set of facets
of $\Phi_{\wA}^{\cF L_r , d}$ which contain $n+1$ is $\{ \sigma \cup
\{ n+1 \}: \sigma \in \Phi_A^{F ,d-1} \}$.
\end{lemma}

{\it Proof}\/. Take $\eta ' \subseteq \{ 1,\ldots ,n+1\}$ such that
$ n+1 \in \eta '$ and set $\sigma=\eta ' \setminus \{ n+1\}$. By
Definition \ref{umbrella}, $\eta ' \in \Phi_{\wA}^{\cF L_r , d}$ if
and only if there is a (unique) vector $\widetilde{c}=(c,c_{d+1})\in
\QQ^{d}\times \QQ$ (which depends on $r$) such that $\langle
\widetilde{c} , \wa_i \rangle = 1$ if $i\in \sigma $, $\langle
\widetilde{c} , -\frac{1}{r} \wa_{n+1} \rangle = 1$ and $\langle
\widetilde{c} , \wa_i \rangle < 1$ if $i\notin \eta '$. The equation
for $i=n+1$ is equivalent to $c_{d+1}= -r$ and the other equalities
and inequalities can then be written as $\langle c , a_i \rangle - r
w_i = 1$ for $i\in \sigma$ and $\langle c , a_i \rangle - r w_i < 1$
for $i\notin  \eta '$. Thus, denoting $c_0 $ the limit of the vector
$c$ when $r$ tends to zero, we have that $\langle c_0 , a_i \rangle
=1$ for $i\in \sigma$ and $\langle c_0 , a_i \rangle \leq 1$ for
$i\in \{1, \ldots , n \} \setminus \sigma$. This proves that $\sigma
\subseteq \eta =\{ i\in \{1, \ldots , n \} :\; \langle c_0 , a_i
\rangle =1\}\in \Phi_A^{F,d-1}$.

It follows for $h=\frac{1}{r}(c-c_0)$ that $\langle h , a_i \rangle
=w_i$ for $i \in \sigma$ and $\langle h , a_i \rangle < w_i$ for
$i\in \eta \setminus \sigma$. This proves that $\sigma \in
\Phi_{A_{\eta}}^{w_{\eta}, d-1}$. This proves the inclusion
$\subseteq$ in the statement. The other inclusion can be proved in a
similar way, now starting from the existence of vectors $c_0$ and
$h$ corresponding to $\eta \in \Phi_A^{F,d-1}$ and $\sigma \in
\Phi_{A_{\eta}}^{w_{\eta}, d-1} $ as above and considering
$\widetilde{c}=(c,-r)\in \QQ^{d}\times \QQ$ with $c= r h + c_0$.
\QED

\begin{lemma}\label{lemma-umbrellas-2}
For any sufficiently large real number $r>0$ we have that
$$\Phi_{\wA}^{\cF L_r ,d}= \{ \sigma \cup \{ n+1 \}: \sigma \in
\Phi_{A_{\eta'}}^{F,d-1}, \eta ' \in \Phi_A^{w,d-1} \}$$ In
particular, if all the facets in $\Phi_A^{w}$ contain exactly $d$
columns of $A$, then $\Phi_{\wA}^{\cF L_r ,d}= \{ \sigma \cup
\{ n+1 \}: \sigma \in \Phi_A^{w,d-1} \}$.
\end{lemma}

{\it Proof}\/.
Recall that $\eta \in \Phi_{\wA}^{\cF L_r ,d}$ if and only if there is a (unique) vector
$\widetilde{c} \in \QQ^{d}\times \QQ$ (which could depend on $r$) verifying the equalities
$\langle \widetilde{c} , \wa_i \rangle = 1$ for $i\in \eta $ and the inequalities
$\langle \widetilde{c} , \wa_i \rangle < 1$ for $i\notin \eta $.

Assume to the contrary that there is a facet $\eta \in
\Phi_{\wA}^{\cF L_r ,d}$ such that $n+1\notin \eta$ for all $r> r_0$
and $r_0 $ big enough. Since $\dim \eta = d$ and $n+1 \notin \eta $
we have that the corresponding $\widetilde{c}$ is independent of
$r>r_0$ and $\langle \widetilde{c} , -\frac{1}{r} \wa_{n+1} \rangle
> 1$ which is  a contradiction. Thus, any facet of $\Phi_{\wA}^{\cF
L_r ,d}$ contains $n+1$.

If we write $\eta = \sigma \cup \{n+1\}$ for $\sigma \subseteq
\{1,\ldots ,n\}$ then $\eta  \in \Phi_{\wA}^{\cF L_r ,d}$ if and
only if there is a unique vector $\widetilde{c}=(c,c_{d+1})\in
\QQ^{d}\times \QQ$ as above (now depending on $r$). In particular,
$c_{d+1}=-r$ and using the limit $c_{\infty}$ of the vector
$\frac{1}{r}c$ when $r$ tends to infinity we get that $\langle
c_{\infty}, a_i \rangle = w_i$ for $i\in \sigma$ and $\langle
c_{\infty}, a_i \rangle \leq w_i $ for $ i\notin \sigma$. Thus,
$\sigma \subseteq \eta '=\{i\in\{1,\ldots,n\} : \; \langle
c_{\infty}, a_i \rangle = w_i \} \in \Phi_A^w$. Moreover, we obtain
that $\sigma \in \Phi_{A_{\eta '}}^F$ by using the vector
$h=c-rc_{\infty}$. The other inclusion can be proved in a similar
way, now starting from the existence of vectors $c_{\infty}$ and $h$
corresponding to $\eta ' \in \Phi_A^{w,d-1}$ and $\sigma \in
\Phi_{A_{\eta '}}^{F, d-1} $ as above and considering $c=  h + r
c_{\infty}$. \QED

\begin{proposition}\label{regularity-along-T}
The following conditions are equivalent:
\begin{enumerate}
\item[(a)] $\{\sigma \in \Phi_{A_{\eta}}^{w_{\eta},d-1}: \; \eta
\in \Phi_A^{F,d-1} \}=\{\sigma \in \Phi_{A_{\eta'}}^{F,d-1}, \eta
'\in \Phi_A^{w,d-1} \}$.
\item[(b)] $\Phi_{\wA}^{\cF L_r}$ is constant for all $r>0$.
\item[(c)] $M_{A , w , \alpha }(\beta )$ does not have slopes along $T$.
\end{enumerate}

\end{proposition}
{\it Proof}\/.
Since any umbrella is determined by its facets, Lemma
\ref{lemma-umbrellas-1} and Lemma \ref{lemma-umbrellas-2} prove that
(a) is equivalent to (b). Corollary \ref{FourierSW} finishes the
proof.
\QED

\begin{remark}
Note that condition (a) in Proposition \ref{regularity-along-T}
implies that there is a common refinement (the one given by
considering the faces of  $\Phi_{\wA}^{\cF L_r}$ not containing
$n+1$) of the polyhedral complex subdivisions induced by the
umbrellas $\Phi_A^F$ and $\Phi_A^{w}$. In particular, when $w$ is
generic, condition (a) means that $w$ induces a regular
triangulation of $A$ that refines the polyhedral complex subdivision
induced by the weight vector $(1,\ldots ,1)$. In other words, $w$ is
a perturbation of $(1,\ldots ,1)$ for the matrix $A$. For example,
this condition is satisfied if either the row span of the matrix $A$
contains the vector $(1, \ldots ,1)$ or $w$ is the sum of a vector
in the row span of the matrix $A$ and a vector $(k , \ldots , k)$
with $k\geq 0$.
\end{remark}

\begin{remark}
Even when $M_{\wA}(\widetilde{\beta})$ is regular holonomic and
$M_{A,w ,\alpha }(\beta )$ has no slopes along $T$, the latter
system is irregular if $M_A (\beta )$ is (see Example
\ref{curious-example}).
\end{remark}

\begin{example}\label{curious-example}
Take $A=(1 \; 2)$ and $\beta \in \CC$. The system
$M_A (\beta)$ is irregular along $Y=\{x_2 = 0\}$ with unique
slope $s=r+1=2$.

If we choose $w =(1,1)$ and consider the matrix $$\wA=\left(\begin{array}{lll}
1 & 2 & 0\\
1 & 1 & 1
\end{array}
\right)$$ we have that the hypergeometric system $M_{\wA} (\beta , \alpha -1 )$ is regular holonomic for all $\beta\in\CC^d$ and $\alpha \in \CC$
by a well known result of Hotta \cite{hotta}.
However, the modified system $M_{A,w ,\alpha}(\beta )$ has a slope $s=r+1=2$ along $T'$ because
$M_{A,-w ,-\alpha}(\beta )$ has the slope $s=2$ along $T$ (see Corollary \ref{FourierSW2} and Remark \ref{Remark-slopes-modify-infinity}).
\end{example}

\begin{remark}
Let $\varphi$ be the map (\ref{eq:gr-deformation}) defined in the
introduction. Since the $\cD$-modules $\cM_{\widetilde{A}({\bf
0})}(\beta , - \alpha )$ and $\cM_{A,w,\alpha}(\beta)$ are
isomorphic when restricted to $X^*=\CC^n\times\CC^*$, the slopes  of
both modules along any coordinate subspace $Z$ not contained in $T$
coincide in $Z\setminus T $. Moreover, the map $\varphi^*$ also
induces an isomorphism for their spaces of Gevrey solutions along
$Z$.
\end{remark}

\subsection{Holomorphic solutions of a modified hypergeometric system}

We study convergent and formal power series solutions of the
modified $A$-hypergeometric module $\cM_{A,w,\alpha}(\beta)$. As was said
before, the map $\varphi$ (\ref{eq:gr-deformation}) induces an
isomorphism between the $\cD$-modules $\cM_{\widetilde{A}({\bf
0})}(\beta , - \alpha )$ and $\cM_{A,w,\alpha}(\beta)$  when
restricted to $X^*$, and also an isomorphism between their
corresponding spaces of holomorphic solutions. More precisely, for
any germ of a holomorphic function $f(x,t)$ at a point $(x_0,t_0)$
in $X^*$, the function $f(x,t)$ is a solution of
$\cM_{A,w,\alpha}(\beta)$ if and only if
$\varphi^*(f)(y,s)=f(s^{-w_1}y_1,\ldots,s^{-w_n}y_n)$ is a germ of a
solution of $\cM_{\widetilde{A}({\bf 0})}(\beta , -\alpha)$ at the
point $(y_0,s_0)\in Y^*$ such that $\varphi(y_0,s_0)=(x_0,t_0)$. We
can rewrite this as follows: the morphism
$$\varphi^* : \cH om_{\cD_{X^*}}(\cM_{A,w,\alpha}(\beta)_{|X^*},\cO_{X^*})
\stackrel{\sim}{\longrightarrow} \cH
om_{\cD_{Y^*}}(\cM_{\widetilde{A}({\bf 0})}(\beta , - \alpha
)_{|Y^*},\cO_{Y^*})$$ is an isomorphism of sheaves of vector spaces.
As a consequence the holonomic ranks of both modules coincide
$$\operatorname{rank}(H_{A,w,\alpha}(\beta)) = \operatorname{rank}(H_{\widetilde{A}({\bf 0})}(\beta , - \alpha ))$$ and
this last rank equals the one of $H_A(\beta)$ for any $w$, see \cite[Theorem 1]{takayama2009}.

Recall that if $\beta$ is generic then $\operatorname{rank}(H_A(\beta))=\operatorname{vol}(A)$,
where $\operatorname{vol}(A)$ is the normalized volume of the matrix $A$ \cite{gkz1989,Adolphson, SST, MMW}, while
in general $\operatorname{rank}(H_A(\beta))\geq \operatorname{vol}(A)$ \cite{SST, MMW}.

\subsection{Gevrey solutions of a modified hypergeometric system}\label{subsection:Gevrey-solutions}

We describe the solutions of $\cM_{A,w,\alpha}(\beta)$ in the
space $\cO_{\widehat{X|T}}$ of formal power series with respect to
$T=\{t=0\}\subset X =\CC^{n+1}$. More generally, we also describe
the solutions of $\cM_{A,w,\alpha}(\beta)$ in the space
$\sum_{\gamma \in \Lambda} t^\gamma \cO_{\widehat{X|T}}$ for any
finite set $\Lambda \subseteq \CC$ (see Theorem \ref{theorem-dimension-1}).

We will use notations in \cite{takayama2009}. Let $\tau$ be the weight vector $({\bf
0},-1,{\bf 0},1) \in \ZZ^{2n+2}$ inducing the Malgrange-Kashiwara $V$--filtration along
$T\equiv (t=0)$ on the ring $D_{n+1}$.

Let $\tilde{I}_{{\widetilde A}(w)}\subseteq \CC [\partial, t]:=\CC
[\partial_1 ,\ldots ,\partial_n , t]$ be the toric ideal associated
with $\wA$, i.e., the binomial ideal generated by the operators in
(\ref{eq:mtoric}).
\begin{lemma}\label{inn_tau-inn_w} For all $w\in\ZZ^n$ we have $\inn_{({\bf 0},-1)}(\widetilde{I}_{{\widetilde A}(w)})= \CC[\partial ,t] \inn_w(I_A)$.
\end{lemma}

{\it Proof}\/. Recall that
$$\widetilde{I}_{{\widetilde A}(w)} = \langle \partial^{u_+} - \partial^{u_-}\, \vert\, Au=0, w\cdot u=0\rangle + \langle \partial^{u_+} - t^{w\cdot u}\partial^{u_-}\, \vert \, Au=0, w\cdot u>0\rangle$$
and we can write $$I_A =\langle \partial^{u_+} - \partial^{u_-}\, \vert\, Au=0, w\cdot u=0\rangle + \langle \partial^{u_+} -\partial^{u_-}\, \vert \, Au=0, w\cdot u>0\rangle .$$
Notice that $\inn_{({\bf 0},-1)}(\partial^{u_+} - t^{w\cdot u}\partial^{u_-}) = \partial^{u_+} = \inn_w(\partial^{u_+} - \partial^{u_-})$ if $Au=0$ and $w\cdot u = w\cdot u_+ -w\cdot u_->0$ and that
$\inn_{({\bf 0},-1)}(\partial^{u_+} - \partial^{u_-}) = \partial^{u_+}- \partial^{u_-} = \inn_w(\partial^{u_+} - \partial^{u_-})$ if $Au=0$ and $w\cdot u = 0$. The conclusion follows by a straightforward Groebner
basis argument because a Groebner basis of $\widetilde{I}_{{\widetilde A}(w)}$ with respect to $({\bf 0},-1)$ (resp. of $I_A$ with respect to $w$) is given by a set of binomials with the same form as the ones defining the ideal.
\QED

Recall that the {\em indicial polynomial} (also called $b$-{\em function}) of
$H_{A ,w}(\beta )$ along $T$ is the polynomial $b(s) \in \CC [s]$ such
that $b(\theta_t )$ is the monic generator of $\inn_\tau(H_{A,w}
(\beta))\cap \CC [\theta_t ]$ where $\theta_t = t \partial_t$.
Moreover, we have by \cite[Th. 3]{takayama2009} that for $\beta$ and
$w$ generic, the indicial polynomial of $H_{A,w}(\beta)$ along $T$
is \begin{equation} b(s)=\prod_{(\partial^{\bf k},\sigma) \in
{\mathcal T}(M)} (s - w \beta^{(\partial^{\bf
k},\sigma)})\label{generic-b-function}\end{equation} where
$M=\inn_{w}(I_A)$, ${\mathcal T}(M)$ is the set of {\em top-dimensional
standard pairs} of $M$ (see \cite[Sec. 3.2]{SST}) and $v=\beta^{(\partial^{\bf k},\sigma)}$ is
the vector defined as $v_i = k_i \in \NN$ for $i\notin \sigma$ and
$A v=\beta$, which is also an exponent of
$H_A(\beta)$ with respect to $w$ (see \cite[Lemma 4.1.3]{SST}).

\begin{definition}\label{perturbation}
We say that a (generic) vector $\widetilde{w}\in \QQ^n$ is a
(generic) perturbation of $w \in \ZZ^n$, with respect to $A$,  if
there exists $w' \in \QQ^n$ such that $\inn_{\widetilde{w}}(I_A) =
\inn_{w'} (\inn_{w}(I_A ))$.\end{definition}

\begin{remark}\label{remark-number-of-exponents}
If $\widetilde{w}$ is generic then $\inn_{\widetilde{w}}(I_A )$ is
a monomial ideal and it is well known that its degree equals the
cardinality of its set of top-dimensional standard pairs ${\mathcal
T}(\inn_{\widetilde{w}} (I_A ) )$. Moreover, for very generic $\beta
\in \CC^d$ there are exactly $\deg (\inn_{\widetilde{w}} (I_A ))$
many exponents of $H_A (\beta )$ with respect to $\widetilde{w}$;
see \cite[Sec. 3.4]{SST} and \cite[Prop. 4.10]{DMM}.
\end{remark}

\begin{lemma}\label{formal-solution-w-perturbation}
Let $\beta \in \CC^d$ be very generic and $w \in
\ZZ^n$. There is a generic perturbation $\widetilde{w}\in \QQ^n$ of
$w$ such that for any exponent $v \in \CC^n$ of $H_A (\beta )$ with
respect to $\widetilde{w}$ the series $\psi_v (x,t)= t^{-\alpha} \phi_v (t^w x)$, for $t^w x=
(t^{w_1 } x_1 , \ldots , t^{w_n } x_n)$, is a solution of
$\cM_{A,w,\alpha}(\beta)$ of the form $\psi_v (x,t)= \sum_{m\geq 0} f_m (x)
t^{\gamma + m }\in t^{\gamma}\cO_{\widehat{X\vert T},(p,0)}$, with $\gamma = w v -\alpha$ and $f_0(x) \neq 0$ for some $p\in \CC^n$.
\end{lemma}

{\it Proof.}\/  Since $\beta$ is very generic, for any generic
$\widetilde{w}$ an exponent $v$ of $H_A (\beta )$ with respect to
$\widetilde{w}$ can be written as $v=\beta^{(\partial^{\bf
k},\sigma)}$ where $(\partial^{\bf k},\sigma)$ is a top-dimensional
standard pair of $\inn_{\widetilde{w}}(I_A )$, see \cite[Sec.
3.4]{SST}. In particular, $\sigma \in \Phi_A^{\widetilde{w},d-1}$,
${\bf k} = (k_i)_{i\not\in \sigma} \in \NN^{n-d}$, $A v=\beta$ and
$v_i = k_i \in \NN$ for all $i\notin \sigma$.

The series $\phi_v (x)$ is either a holomorphic solution or a Gevrey solution of
$\cM_A (\beta )$ along a coordinate subspace $Z\subseteq \CC^n$ at any point $p$
in a non empty relatively open set $U_\sigma$ in $Z$ (see Theorem
\ref{Theorem-Gevrey-solutions-hypergeometric} and \cite[Th. 3.11]{Fer} for the details) and
since $\beta$ is very generic
we have that $N_v = (-B_{\sigma} {\bf k} + \NN B_{\sigma})\cap \ZZ^n$.

The  expression  $f(x,t):= t^{-w v} \phi_v (t^{w_1 } x_1 , \ldots , t^{w_n }x_n )=t^{\alpha -wv}\psi_v (x,t)$
(resp. $\psi_v(x,t)$) formally satisfies the equations defining $\cM_{A ,w ,w v}(\beta)$
(resp. $\cM_{A ,w ,\alpha}(\beta)$). We will prove that we can write $f(x ,t)=\sum_{m\geq 0} f_m(x) t^{m}$ and that it is a Gevrey series along $T \subset X$. Recalling the expression of $\phi_v$ (\ref{phiv}) it is enough to prove that for all $u \in N_v \setminus \{0\}$ we have $w u\in
\NN $ and that the coefficient of $t^m$ in $f$, i. e. $$f_m
(x)=\sum_{u \in N_v , w u =m} \dfrac{[v]_{u_{-}}}{[v+u]_{u_{+}}}
x^{v+u},$$ is a convergent series in an open neighborhood of some $p\in \CC^n$, both the neighborhood and $p$ independent of $m$.

We can take a generic perturbation $\widetilde{w}\in \QQ^n$ of $w$
of the form $\widetilde{w}=w + \epsilon \widetilde{e}$ with $\widetilde{e}=(1, \ldots ,1 ) + \epsilon'w' $ for $\epsilon >0$ and $\epsilon '>0$ small
enough and $w'\in \QQ^n$ is generic.

Take any $u\in N_v \setminus \{0\}$ and let us prove that $w u\geq 0$.
Since $v$ is an exponent of $H_A (\beta )$ with respect to
$\widetilde{w}$ we have by \cite[(3.30)]{SST} that $\widetilde{w} u
> 0$. Hence, since last inequality holds for $\epsilon >0$ and $\epsilon
' >0$ small enough we have that $w u \geq 0$. Thus we have that $w u \geq 0$ for
all $u\in N_v$. Notice that when $w u
> 0$ for all $u\in N_v \setminus \{0\}$ then the set $\{u \in N_v ,
w u =m \}$ is finite and hence $f_m (x)$ is clearly convergent. In
general, $\{u \in N_v ,
w u =m \}$ is not finite, but we will see that $f_m (x)$ is still convergent at some point $p\in \CC^n$.
Since $N_v = (-B_{\sigma} {\bf k} + \NN B_{\sigma})\cap \ZZ^n$ the set
$\{u \in N_v ,w u =m \}$ is a finite union of shifted copies of the form
$N(i)=u(i) + (\sum_{j\notin \sigma; w b_j=0} \NN b_j )\cap \ZZ^n$ where $\{ b_j :\; j \notin \sigma\}$ is the set of
columns of $B_{\sigma}$ and $u (i) \in N_v$ satisfies $w u(i) =m$. The series $f_m (x)$ is
convergent if and only if all the series $g_{i,m}(x)=\sum_{u\in N(i)}  \dfrac{[v]_{u_{-}}}{[v+u]_{u_{+}}}
x^{v+u}$ are convergent. Since $\beta$ is very generic, $\operatorname{nsupp} (v)=\emptyset$, and
the convergence of each series $g_{i,m}$ is equivalent to the convergence of the series
$x^{v+u(i)}\sum_{u\in -u(i) + N(i)}\dfrac{|u_{-}|!}{|u_{+}|!} x^{u}$. Thus, it is enough to see that
for any column $u$ of $B_{\sigma}$ such that $w u=0$ we have that $|u|=|u_{+}|-|u_{-}|\geq 0$. Notice that $w u
=0$ implies $0< \widetilde{w} u =  \epsilon (|u| + \epsilon ' w ' u)$ and so $|u| + \epsilon ' w ' u >0$.
Hence, since this holds for $\epsilon '  >0$ small
enough, we have that $|u|\geq 0$. We have proved that $f$ is a
formal solution of $M_{A,w,wv}(\beta )$ along $T$ and it is clear
that $f_0 (x)=x^v + \cdots \neq 0 $.
From the expression of the $g_{i,m}$ and \cite[Th. 3.11]{Fer} any $f_m(x)$ is convergent at any point in $\{x\in \CC^n\, \vert\, 0\not=\prod_{i\in \sigma} x_i, |x_j| <R|x_\sigma^{A_\sigma^{-1}a_j}| {\mbox{ for }} j\not\in \sigma {\mbox{ and }}  |wb_j|=0\}$ for some $R>0$.
\QED

Let $\widetilde{w}\in \QQ^n$ be a generic perturbation of $w \in \ZZ^n$ as in the proof of Lemma \ref{formal-solution-w-perturbation}.

\begin{lemma}\label{f_m} If $f(x,t)=\sum_{m\geq 0} f_m(x) t^{\gamma+m}\in t^{\gamma}\cO_{\widehat{X\vert T},p}$ is a
solution of  $\cM_{A,w,\alpha}(\beta)$ for some $\gamma \in \CC$, $p\in T$, with $f_0(x) \neq  0$, then:
\begin{enumerate}
 \item[(a)] $t^{\alpha}f(x,t)$ is a solution of $\cM_{A,w}(\beta)$.
 \item[(b)] For all $m\geq 0$, $f_m(x)$ is a holomorphic solution of $\cM_{A_w}(\beta, \alpha+\gamma +m)$, where this last module is the hypergeometric system associated with  the matrix $A_w$ and the parameter $(\beta, \alpha+\gamma +m)$.
 \item[(c)] $b(\alpha + \gamma )=0 $, where $b(s)$ is the indicial polynomial of $H_{A,w}(\beta)$ along $T$.
 \item[(d)] If $\beta$ is very generic then $\alpha + \gamma = w v$ for some exponent $v$ of $H_A (\beta )$ with respect to $\widetilde{w}$.
\end{enumerate}

\end{lemma}

{\it Proof.\/}
The proof of (a) and (b) are straightforward. Let us prove
(c). By (a) and using \cite[Theorem
2.5.5]{SST} we have that $\inn_{(0,1)}(t^{\alpha}f(x,t))=f_0 (x) t^{\alpha
+\gamma }$ is a solution of $\inn_\tau (H_{A,w} (\beta))$.

Recall by definition of $b(s)$ that $\langle b(\theta_t ) \rangle
=\inn_\tau (H_{A,w} (\beta))\cap \CC [\theta_t ]$ where $\theta_t =
t\partial_t$. Thus, the differential operator $b(\theta_t )$
annihilates $\inn_{(0,1)}(t^{\alpha}f(x,t))=f_0 (x) t^{\alpha
+\gamma }$. This implies, using for example \cite[Lemma 1.3.2]{SST},
that $0= b(\theta_t )(f_0 (x) t^{\alpha +\gamma})=b(\alpha
+\gamma)f_0 (x) t^{\alpha +\gamma} $ and this implies that $b(\alpha
+\gamma)=0$.

Let us prove (d). By (b) we have that $f_0 (x)$ is a holomorphic
solution of $H_{A_w}(\beta, \alpha+\gamma )$ and thus it can be
written as a Nilsson series at the origin with respect to a vector
$\widetilde{e}$ that is a perturbation of $e=(1 ,\ldots ,1)$ (see
\cite{SST}, \cite{ohara-takayama}, \cite{DMM}) and, in particular,
it makes sense to consider the initial form of $f_0 (x)$ with
respect to $\widetilde{e}$. On the other hand, using Lemma
\ref{inn_tau-inn_w}, we have that $\inn_w (I_A )+ \langle A \theta
-\beta \rangle \subseteq \inn_\tau (H_{A,w} (\beta))$ annihilates
$f_0 (x)$. Since $\beta$ is very generic, $\inn_w (I_A )+ \langle A
\theta -\beta \rangle = \inn_{(-w,w)}H_A (\beta)$ \cite[Th. 3.1.3]{SST}. This implies that
$f_0 (x)$ is a solution of $\inn_{(-w,w)}H_A (\beta)$ and hence,
$\inn_{\widetilde{e}}(f_0(x))$ is a solution of
$\inn_{(-\widetilde{w},\widetilde{w})}H_A (\beta)$ for
$\widetilde{w}=w + \epsilon \widetilde{e}$. Thus, since $\beta$ is
very generic $\inn_{\widetilde{e}}(f_0(x))=c x^v$ for $c\in \CC$ and
$v$ an exponent of $H_A (\beta)$ with respect to $\widetilde{w}$.
Hence, using (b), we also have that $w v= \alpha + \gamma$.
\QED

\begin{remark}
Although we assume in this paper that $\wA$ is pointed it turns out that in this section this fact is only used in the proof of (d) in Lemma \ref{f_m}.
However, let us notice that if $\wA$ is not pointed and $w$ is generic then $\inn_w (I_A )=\CC [\partial]$, $\inn_\tau (H_{A,w} (\beta))=D$ and so $b(s)=1$.
Thus, by (c) in Lemma \ref{f_m} the modified system $\cM_{A,w,\alpha}(\beta)$ does not have any solution in $t^{\gamma}\cO_{\widehat{X\vert T},p}$ for all $\gamma \in \CC$ and $p\in T$.
       \end{remark}

\begin{remark}\label{degree-w-initial-toric}
Since $\widetilde{w}=w + \epsilon \widetilde{e}$ with
$\widetilde{e}=(1, \ldots ,1 ) + \epsilon ' w ' $ for sufficiently small $\epsilon >0$
and $\epsilon '>0$ we have that $\inn_{\widetilde{w}}
(I_A )=\inn_{w'}(\inn_e(\inn_w (I_A )))$ for $e=(1 ,\ldots , 1)$. In
particular, $\inn_{\widetilde{w}} (I_A )$ and $\inn_w (I_A )$ have
the same degree.
\end{remark}

Let us denote by $\dim_\CC(M,{\cal F})_p$ the dimension of the space of ${\cal F}$-solutions of a $D$-module $M$ at a point $p$.

\begin{theorem}\label{theorem-dimension-1} Assume $\beta\in \CC^d$ is very generic, $w \in \ZZ^n$ and $\alpha\in \CC$.
Then $\dim_\CC(\cM_{A,w,\alpha}(\beta), \cO_{\widehat{X\vert T}})_p
=0$ if $w v -\alpha \notin \NN$ for all the exponents $v$ of $H_A
(\beta)$ with respect to $\widetilde{w}$. We also have that
$$\dim_\CC(\cM_{A,w ,\alpha }(\beta), \sum_{b(\alpha + \gamma )=0} t^{\gamma} \cO_{\widehat{X\vert T}})_{p} =\deg (\operatorname{in}_{w}(I_A)) $$

In particular, if $w$ is generic we also have $$\dim_\CC(\cM_{A,w,wv-m}(\beta), \cO_{\widehat{X\vert T}})_p =1$$ for
all generic $p\in T$, $m \in \NN$ and any exponent $v$ of $H_A (\beta )$ with respect to $w$.
\end{theorem}

{\it Proof}\/.
The first statement follows from Lemma \ref{f_m}, (d).

The inequality $\dim_\CC(\cM_{A,w ,\alpha }(\beta), \sum_{b(\alpha + \gamma )=0} t^{\gamma} \cO_{\widehat{X\vert T}})_p \geq \deg (\operatorname{in}_{w}(I_A))$
follows from Lemma \ref{formal-solution-w-perturbation}, Remark \ref{remark-number-of-exponents} and Remark \ref{degree-w-initial-toric}. On the other hand,
since the differential operators defining $\cM_{A,w,\alpha}(\beta)$ belong to the Weyl Algebra we have that any solution $f\in \sum_{b(\alpha + \gamma )=0} t^{\gamma} \cO_{\widehat{X\vert T}}$
of $\cM_{A,w,\alpha}(\beta)$ decomposes as a finite sum of solutions, each of them in a space $t^{\gamma} \cO_{\widehat{X\vert T},p}$. Recall by the proof of Lemma \ref{f_m}
that any solution $f \in t^{\gamma} \cO_{\widehat{X\vert T},p}$ of $\cM_{A,w,\alpha}(\beta)$ verifies that $\inn_{(\widetilde{w}, 0)} (\inn_{({\bf 0},1)} (f))=c x^v t^{w v -\alpha}$
for an exponent of $H_A (\beta)$ with respect to $\widetilde{w}$. This last fact, Remark \ref{remark-number-of-exponents}, Remark \ref{degree-w-initial-toric} and a slightly modified version of \cite[Proposition 2.5.7]{SST}
prove that $\dim_\CC(\cM_{A,w ,\alpha }(\beta), \sum_{b(\alpha + \gamma )=0} t^{\gamma} \cO_{\widehat{X\vert T}})_p \leq \deg (\operatorname{in}_{w}(I_A))$.

Finally, the last statement follows from the second statement and from the
fact that if $\beta$ is very generic, $w\in
\ZZ^n$ is generic and $v$ and $v'$ are two different
exponents of $H_A (\beta )$ with respect to $w$, then $w (v-v')\notin \ZZ$.
\QED

\begin{remark}
\label{remark3} Let $\psi_v(x,t)$ be the series constructed in Lemma
\ref{formal-solution-w-perturbation} and used in Theorem
\ref{theorem-dimension-1}. If $w$ is in the row span of $A$ then $f
(x,t)=t^{\alpha -wv} \psi_v$ does not depend on $t$ and thus it is a
convergent series. If $w$ is not in the row span of $A$ then $f
(x,t)$ is Gevrey along $T$ with index $s=r+1$ where
$$r=\max \{-\dfrac{|u|}{w u} :\; \NN u \subseteq N_v, w u> 0 \}$$ where $|u|=\sum_i u_i$. On the other hand, as mentioned in the proof of Lemma \ref{formal-solution-w-perturbation} since $\beta$ is very generic and $v$ is an exponent of $H_A (\beta )$ with respect to $\widetilde{w}$ then $v$ is associated with a simplex
$\sigma \in \Phi_A^{\widetilde{w}, d-1}$
and there is a basis $\{b_i :\; i\notin \sigma \}$ of the kernel of $A$ such that for all $i\notin \sigma$, $(b_i)_j =0$ for all $j\notin \sigma \cup \{ i\}$ and $(b_i)_i =1$.
The set $\{b_i :\; i\notin \sigma \}$ is the set of columns of $B_{\sigma}$ (if we reorder the variables so
that $\sigma =\{1 ,\ldots ,d\}$) and in this case we have that $N_v = (-B_{\sigma} {\bf k} + \NN B_{\sigma})\cap \ZZ^n$.
Thus, more explicitly, $r=\max \{-|b_i |/(w b_i) :\; i \notin \sigma, w b_i >0 \}$ where $\{b_i : \; i\notin \sigma \}$ is the set of columns of $B_{\sigma}$,  $|b_i |=1-|A_{\sigma}^{-1}a_i|$ and $w b_i =
w_i -w_{\sigma} A_{\sigma}^{-1} a_i >0$. The proof of this formula is technical and follows from standard estimates on Gamma functions similar to the ones used in \cite{Fer} to compute
the index of Gevrey solutions for hypergeometric systems. In particular, if $w$ is a perturbation of $(1,\ldots ,1)$ then $r$ is close to $-1$ and if $A$ is homogeneous then $r=0$ because $|u|=0$ for
any $u\in N_v$ and hence in both cases the series is convergent. \end{remark}

\subsection{Gevrey solutions modulo convergent series}\label{Gevrey-solutions-modulo-convergent-series}

By Theorem \ref{theorem-dimension-1}, if both $\alpha
\in \CC$ and $\beta \in \CC^d$ are very generic then $\cM_{A ,w
,\alpha} (\beta )$ does not have any nonzero solution in
$\cO_{\widehat{X\vert T},p}$ for all $p\in T$. This is in contrast
with the case of the irregularity of hypergeometric systems along
coordinate hyperplanes, where for any slope $s=r+1$ of $\cM_A (\beta
)$ along $Y=\{x_n = 0\}$ and for very generic $\beta
\in \CC^d$ one can construct a formal solution
$\phi \in \cO_{\widehat{X\vert Y},p}$ of $\cM_A(\beta )$ along $Y$,
such that $\phi$ has Gevrey index equal to the slope (see \cite{Fer}
and Theorem \ref{hypergeometric-Gevrey-solutions}).

However, by the comparison theorem for the slopes
\cite{Laurent-Mebkhout-pa-pa} and the perversity of the irregularity
complex of a holonomic $\cD$-module along a smooth hypersurface
\cite{Mebkhout-positivite}, one knows that for each slope
$s=r+1$ of $\cM_{A, w ,\alpha }(\beta )$ along $T$ at a generic
$p\in T$ there must exist a formal series $\phi \in
\cO_{\widehat{X\vert T},p}$ with Gevrey index $s=r+1$ such that
$P(\phi )$ is convergent at $p$ for all $P\in H_{A ,w ,\alpha}
(\beta )$.

The purpose of this section is to describe Gevrey solutions modulo
convergent series of the modified system $\cM_{A ,w ,\alpha
}(\beta)$ along $T$ when $\alpha \in \CC$ and $\beta \in \CC^d$ are
very generic. To this end, we use the construction of the
Gevrey solutions at infinity of the hypergeometric system
$M_{\wA}(\beta , \alpha -1 )$ as performed in Theorem
\ref{hypergeometric-Gevrey-solutions-infinity}.

Take $X'=\CC^{n+1}$ with coordinates $(x_1 ,\ldots ,x_n , z)$ and
$z=1/t$ so that $X\cap X'=\CC^n \times \CC^{\ast}$. Denote
$T'=\{t=\infty\}=\{z=0\}\subseteq X'$. We can consider for any
$\gamma \in \CC$ the $\CC$--linear map
$$\begin{array}{rcl}
   \Upsilon_{\gamma} : t^{\gamma }\cO_{\widehat{X|T},p} & \longrightarrow & t^{-1-\gamma} \cO_{\widehat{X'|T'},p'}\\
 f=\sum_{m\geq 0} f_{m}(x) t^{\gamma +m} & \longmapsto & \Upsilon_{\gamma} (f)=\sum_{m\geq 0} f_{m}(x) [-\gamma - 1]_{m} t^{-1-\gamma -m}
  \end{array}
$$ where $p=(p_1 ,\ldots , p_n ,0)\in T$ and $p'=(p_1 ,\ldots , p_n , \infty  )\in T'$.

\begin{remark}\label{upsilon-morphism}
Notice that $\Upsilon_{\gamma}$ is an isomorphism if and only if
$\gamma \notin \ZZ_{<0}$. In such a case we also have that
$\Upsilon_{\gamma}(t^{\gamma }\cO_{\widehat{X|T}}(s-1))=
t^{-1-\gamma} \cO_{\widehat{X'|T'}}(s)$ for all $s$. It is also
clear that $\Upsilon_{0}(\sum_{m\geq 0} f_{m}(x) t^{k +m}) =[-1]_k
\Upsilon_{k}(\sum_{m\geq 0} f_{m}(x) t^{k +m})$ for all $k \in \NN$.
\end{remark}

\begin{theorem}\label{formal-modulo-convergent}
Assume $\alpha \in \CC$ and $\beta \in \CC^d$ to be very generic. If
$s=r+1>1$ is a slope of $\cM_{A ,w ,\alpha}(\beta )$ along $T$ then
we can construct
$\sum_{\tau}\operatorname{vol}(\operatorname{conv}(0, \wa_i : \;
i\in \tau ))$ Gevrey series that are linearly independent solutions
of $\cM_{A ,w ,\alpha}(\beta )$ modulo convergent series and whose
Gevrey index is equal to $s=r+1$. Here $\tau$ runs over all the
facets of $\Delta_{A_w}$ such that $-\frac{1}{r}\wa_{n+1} \in
H_{\tau}$ and $0\notin \tau$. Moreover, the classes modulo $
\cO_{\widehat{X|T}}({<}\,\!s)$ of these Gevrey series form  a basis of the
solution space of $\cM_{A,w,\alpha}(\beta )$ in
$(\cO_{\widehat{X|T}} (s)/\cO_{\widehat{X|T}}({<}\,\!s))_p$ for points $p\in T$ in a relatively open set of $T$.\end{theorem}

{\it Proof}\/. The existence of such facets $\tau$ is given by
Corollary \ref{FourierSW2}. Since $-\frac{1}{r}\wa_{n+1} \in
H_{\tau}$ and $-r=2-s'$ for $s'=s+1>2$ we have that $s'>2$ is a
slope of $\cM_{\wA}(\beta , \alpha -1)$ along $T'=\{t=\infty\}$.
Thus, by Theorem \ref{hypergeometric-Gevrey-solutions-infinity} we
can construct $\sum_{\tau}\operatorname{vol}(\operatorname{conv}(0,
\wa_i : \; i\in \tau )$ Gevrey series along $T'=\{t= \infty \}$ with
index $s '$. Moreover, the classes in
$\cO_{\widehat{X'|T'}}(s')/\cO_{\widehat{X'|T'}}({<}\,\!s')$ are
linearly independent solutions of $\cM_{\wA}(\beta , \alpha -1)$.

More precisely, for any $d$-simplex $\sigma \subseteq \tau$ the
series constructed are of the form $\phi_{\widetilde{v}}$ for
$\widetilde{v}=(v,-1-k)$ with $A v=\beta$, $w v -1 -k = \alpha -1$
(i. e. $w v -\alpha = k \in \NN$) and $v_i \in \NN$ for all $i\in
\{1 ,\ldots , n\}\setminus \sigma$.

Using Remark \ref{upsilon-morphism} we can take $\psi_{v} (x,t )$ as
the unique Gevrey series along $T$ with index $s=r+1$ verifying
$\Upsilon_{0} (\psi_{v} (x,t ))=\phi_{\widetilde{v}}$ for
$\widetilde{v}=(v,-1-k)$.

We conclude by Remark \ref{upsilon-morphism} that the
$\sum_{\tau}\operatorname{vol}(\operatorname{conv}(0, \wa_i : \;
i\in \tau )$ series $\psi_{v} (x,t )$ constructed are  Gevrey series
with index $s=s'-1=r+1$ whose classes modulo $\cO_{\widehat{X|T}}({<}\,\!s)$ are
linearly independent. Moreover, it can be checked that they are
solutions of the modified system modulo $\cO_{X|T}$ by using the
fact that their images by the morphism $\Upsilon_{0}$ are solutions
of $\cM_{\wA}(\beta , \alpha -1)$.

Last statement follows from (\ref{FouriermGKZ=GKZ}),
\cite{Laurent-Mebkhout-pa-pa} and \cite{SW}. \QED

\begin{example}
Take $A=(1\; 3 \; 5)$, $w =(0,1,1)$ and $\beta , \alpha \in
\CC$. We have that $\widetilde{I}_{\wA} = \langle \partial_2 - t
\partial_1^3 , \partial_3 - t \partial_1^5 \rangle$ and $H_{A,
w , \alpha }(\beta )= D \widetilde{I}_{\wA} + D \langle x_1
\partial_1 + 3 x_2 \partial_2 + 5 x_3 \partial_3 -\beta , x_2
\partial_2 + x_3 \partial_3  -t \partial_t -\alpha \rangle$. Note here that $\widetilde{I}_{\wA}$ is the binomial ideal generated by the operators in
(\ref{eq:mtoric}). The
unique slope of $M_{A ,w ,\alpha }(\beta )$ along $T$ is
$s=r+1=5$ since $-\dfrac{1}{4}\wa_4$ belongs to the line passing
through $\wa_1 ,\wa_3$ and $\sigma =\{1,3\}$ is a facet of
$\Phi_{\wA}^{\cF L_{r}}$ if and only if $r\geq 4$.

The volume of $\sigma$ is one and following the proofs of Theorem
\ref{formal-modulo-convergent} and Theorem
\ref{hypergeometric-Gevrey-solutions-infinity} we can take $v=(\beta
- 5 \alpha ,0 , \alpha )$ which satisfies the conditions $w v
-\alpha = k = 0 \in \NN$, $v_2 = 0 \in \NN$ and $ A v= \beta$. We
get that the series
$$\psi_v (x ,t )= \sum_{m_2 , m_2 + m_3 \geq 0} \dfrac{[\beta - 5\alpha]_{3 m_2 + 5 m_3}}{[\alpha +m_3]_{m_3} m_2 !} x_1^{\beta -5 \alpha -3m_2 -5 m_3} x_2^{m_2} x_3^{\alpha + m_3} t^{m_2 + m_3} $$
is a Gevrey solution (modulo convergent series) of $M_{A ,w
,\alpha }(\beta )$ along $T$ with index $s=r+1=5$.
\end{example}

\section{Borel transformation and asymptotic expansion}
\label{section:Borel}
We assume that the $\QQ$-row span of the matrix $A$ does not contain
the vector $(1, \ldots, 1)$,
but the one of $A_w$ does,
where $A_w$ is the matrix with columns $\wa_i$, $1\leq i\leq n$ (see Subsection \ref{slopes-mhgs}).
This case holds if and only if the weight vector
$w$ is in the image of ${\bar A}^T$ where ${\bar A}$ is the matrix
$\left(
\begin{array}{ccc}
 a_1 & \cdots & a_n \\
  1  & \cdots & 1 \\
\end{array} \right)
$ and it is not in the row span of $A$. In other words,
the weight vector $w$ lies in the intersection of the set of the
secondary cones of ${\bar A}$; see \cite[Ch. 8]{Sturmfels}.

Solutions of this case can be analyzed by utilizing the Borel transformation
and the Laplace transformation.

We review here some basics of the Borel summation
method which we require in the following (see \cite{B} for the details).
Let us consider the formal expression
\begin{equation}
\label{borel-formal}
f(t) = \sum_{\ell = 0}^{\infty} f_{\ell} t^{\ell + \gamma}
\in t^{\gamma}\mathbb{C}[[t]]
\end{equation}
where $f_0\not= 0$ and $\gamma\in \CC$.
Solutions constructed in Theorem \ref{theorem-dimension-1} are of this form.
If its coefficients satisfy
\begin{equation}
\label{Gevrey-estimates}
\big| f_{\ell} \big| \leq C K^\ell \Gamma(1+ (\ell + \gamma)/\kappa) \qquad (\ell = 0, 1, 2, \cdots)
\end{equation}
with some positive constants $C, K, \kappa$, and $\Re \gamma > -\kappa$ (in \eqref{borel-expnent}
this last condition will be relaxed),
then the formal Borel transform
(with index $\kappa$) defined by
$$
\hat{\mathcal{B}}_{\kappa}[f](\tau)
:= \sum_{\ell = 0}^{\infty}
\frac{f_{\ell}}{\Gamma(1 + (\ell + \gamma) /\kappa)}
\tau^{\ell + \gamma}
$$
is the product of $\tau^\gamma$ and a convergent power series at $\tau = 0$.
In addition to \eqref{Gevrey-estimates}, if
\begin{itemize}
\item[(i)]
the function
$\hat{\mathcal{B}}_{\kappa}[f]$ can be analytically continued to a sector
$$
S(\theta, \delta)
:=\{r e^{i\theta'};\, \big|\theta' - \theta\big| < \delta/2,
\, r > 0
\}
$$
of infinite radius in a direction $\theta \in \mathbb{R}$
with an opening angle $\delta > 0$, and

\item[(ii)]
the analytic continuation of
$\hat{\mathcal{B}}_{\kappa}[f]$ satisfies the growth estimate
\begin{equation}
\label{borel-summability-ii}
\Big|\hat{\mathcal{B}}_{\kappa}[f](\tau)\Big| \leq c_1 \exp\big[c_2 \big|\tau\big|^{\kappa}\big]
\end{equation}
in $S(\theta, \delta)$ with some positive constants $c_1, c_2 > 0$,
\end{itemize}
then we say $f$ is $\kappa$-summable in the direction $\theta$,
and define the $\kappa$-sum
(or the Borel sum with index $\kappa$) of $f$ by
the Laplace transformation
\begin{equation}
\label{def-Borel-sum}
\mathcal{S}[f](t)
= \mathcal{L}^{\theta}_{\kappa}\hat{\mathcal{B}}_{\kappa}[f](t)
:=
\int^{e^{i\theta}\cdot \infty}_0 e^{-(\tau/t)^\kappa} \hat{\mathcal{B}}_{\kappa}[f](\tau) d(\tau/t)^{\kappa},
\end{equation}
where
$$
d(\tau/t)^{\kappa} = \dfrac{\kappa\tau^{\kappa-1}}{t^{\kappa}} d\tau.
$$

\begin{remark}\label{remak(a)}  Because of the growth condition (ii) of $\hat{\mathcal{B}}_{\kappa}[f]$,
the Laplace integral \eqref{def-Borel-sum} converges if
$t$ satisfies
\begin{equation}
\label{def-Borel-sum-domain-tmp}
 \Re\left[\left(\frac{\tau}{t}\right)^{\kappa}\right] - c_2 |\tau|^{\kappa} > 0.
\end{equation}
Since
$$
\Re\left[\left(\frac{\tau}{t}\right)^{\kappa}\right] - c_2 |\tau|^{\kappa}
=
 \left|\frac{\tau}{t}\right|^{\kappa}
\big\{
\cos \kappa(\theta - \arg t) - c_2 |t|^{\kappa}\big\}
$$
(note that $\arg \tau = \theta$),
the Laplace integral \eqref{def-Borel-sum} converges in
\begin{equation}
\label{def-Borel-sum-domain}
\{t;\, \cos\big[\kappa( \arg t - \theta)] > c_2 |t|^{\kappa}\}.
\end{equation}
The region \eqref{def-Borel-sum-domain}
has infinitely many connected components.
Here and in what follows we specify one of them by
imposing $|\arg t - \theta| < \pi/(2 \kappa)$.
Since we can vary $\arg \tau$ in \eqref{def-Borel-sum} slightly,
we conclude that ${\mathcal S}[f]$ defines a holomorphic function in
$$
\bigcup_{|\theta' - \theta| < \delta/2}
\{t;\, \cos\big[\kappa(\arg t - \theta')] > c_2 |t|^{\kappa},\,
|\arg t - \theta'| < \pi/(2\kappa)\}.
$$
Therefore we can find $\rho > 0$ and $\varpi > \pi/\kappa$
such that the ${\mathcal S}[f]$ is holomorphic in
$S(\theta, \varpi, \rho)
:=S(\theta, \varpi) \cap \{t;\, 0 < |t| < \rho\}$
(cf. \cite[the first paragraph of \S 2.1]{B}).
\end{remark}

\begin{theorem}\label{old-remark-(b)}
If $f$ is $\kappa$-summable in a direction $\theta$,
then $f$ is a Gevrey asymptotic expansion of its Borel sum $\mathcal{S}[f](t)$:
For any closed subsector $\overline{S}$ of  $S(\theta, \varpi, \rho)$,
there exists $C', K' > 0$ for which the inequality
\begin{equation} \label{eq:asymptotic_expansion}
\left|
t^{-\gamma}
\mathcal{S}[f](t) - \sum_{\ell = 0}^{N-1} f_{\ell} \, t^{\ell}
\right|
\leq C' (K')^N \big|t\big|^N \Gamma(1 + N/\kappa)
\end{equation}
holds in $\overline{S}$ for $N \in \NN$. \end{theorem}
{\it Proof}\/.
Here we give a sketch of the proof.
See \cite[Th. 1]{B} for the details.

It follows from the relation
$$
t^{\ell + \gamma} =
\int^{e^{i\theta}\cdot \infty}_0 e^{-(\tau/t)^\kappa}
\dfrac{\tau^{\ell + \gamma}}{\Gamma(1 + (\ell + \gamma)/\kappa)}
\dfrac{\kappa\tau^{\kappa-1}}{t^{\kappa}} d\tau
\quad
(= \mathcal{L}^{\theta}_{\kappa}\hat{\mathcal{B}}_{\kappa}[t^{\ell + \gamma}])
$$
that the remainder of the expansion becomes
\begin{align}
&
\mathcal{S}[f](t) - \sum_{\ell = 0}^{N-1} f_{\ell} \, t^{\ell + \gamma}\label{old-remark-(b)-key-relation}\\
&
\notag
\qquad
=
\int^{e^{i\theta}\cdot \infty}_0 e^{-(\tau/t)^\kappa}
\left\{
\hat{\mathcal{B}}_{\kappa}[f] (\tau) -
\sum_{\ell = 0}^{N-1} \frac{f_{\ell}}{\Gamma(1 + (\ell + \gamma)/\kappa)} \, \tau^{\ell + \gamma}
\right\}
\dfrac{\kappa\tau^{\kappa-1}}{t^{\kappa}} d\tau.
\end{align}
Since
$$
\frac{1}{\tau^{N + \gamma}}
\left\{
\hat{\mathcal{B}}_{\kappa}[f] (\tau) -
\sum_{\ell = 0}^{N-1} \frac{f_{\ell}}{\Gamma(1 + (\ell + \gamma)/\kappa)} \, \tau^{\ell + \gamma}
\right\}
$$
is holomorphic in some neighborhood of $\overline{S}$ (which includes the origin)
and satisfies the growth condition \eqref{borel-summability-ii} in $\overline{S}$
with appropriate constants $c_1$ and $c_2$,
the righthand side of \eqref{old-remark-(b)-key-relation} can be
estimated by the righthand side of \eqref{eq:asymptotic_expansion}
multiplied by $|t|^{\gamma}$.
\QED

Inequality \eqref{eq:asymptotic_expansion} is also known to be equivalent to conditions (i) and (ii) stated above (\cite[Theorem 1 (p. 23)]{B}).

The Borel summation method may be also applied
to the case when
\begin{equation}
\label{borel-expnent} \gamma \not\in \kappa \ZZ \quad\text{and}\quad \ell + \gamma \not\in -\kappa\NN_{>0}  \quad\text{for}\quad \ell \in \NN.
\end{equation}
We can define the Borel transform $\hat{\mathcal{B}}_{\kappa}[f]$
in the same manner as before.
In this last case, however,
the Laplace integral \eqref{def-Borel-sum} may not converge
at $\tau = 0$.
Therefore we modify the definition of the Borel sum to
\begin{equation}
\label{def-Borel-sum2}
\mathcal{S}[f](t)
=
\mathcal{L}^{\theta}_{\kappa}\hat{\mathcal{B}}_{\kappa}[f](t)
:=
\frac{1}{1 - e^{- 2\pi i \gamma/\kappa}}
\int_{\Gamma_{\kappa\theta}}
e^{-\zeta/t^\kappa}
\hat{\mathcal{B}}_{\kappa}[f](\zeta^{1/\kappa}) \dfrac{d\zeta}{t^{\kappa}}
\end{equation}
with a path of integration $\Gamma_{\kappa\theta}$
which runs from $\infty$ along $\arg \zeta = \kappa\theta -2 \pi$
to some point near the origin,
takes a $2 \pi$ radian turn along a circle with the center at the origin,
and goes back to infinity in the direction $\arg \zeta = \kappa\theta$.
When condition \eqref{borel-expnent} is satisfied and $\Re \gamma > -\kappa$, then
\eqref{def-Borel-sum2} coincides with \eqref{def-Borel-sum}.
The Borel sum \eqref{def-Borel-sum2} also satisfies the same properties as previously defined \eqref{def-Borel-sum}.
\bigskip

We apply the Borel summation method to the Gevrey
solution constructed in Theorem \ref{theorem-dimension-1}.
Our main result of this section is

\begin{theorem}
\label{borel-thm}
We assume that the $\QQ$-row span of the matrix $A$ does not contain
the vector $(1,1, \ldots, 1)$ but that the one of $A_w$ does. We also assume $\beta$ to be very generic.

Let \begin{equation}
\label{formal-solution}
\psi(x, t) = \sum_{\ell = 0}^\infty C_{\ell} (x) t^{\ell + \gamma}
\end{equation}
be one of the formal solutions of the modified hypergeometric system
$H_{A, w}(\beta)$ constructed in Theorem \ref{theorem-dimension-1}
and $r+1$ be the Gevrey index
of $\psi(x, t)$ along $T$. We
also assume $r\gamma\not\in\mathbb{Z}$. Then the formal solution
$\psi(x, t)$ is $1/r$-summable (as a formal power series in $t$) in
all but finitely many directions for each $x \in U$ where
$U$
is a non-empty open set in the $x$-space $\CC^n$. Furthermore its Borel sum
determines a solution of the modified hypergeometric system $H_{A,
w}(\beta)$.
\end{theorem}

\begin{remark}\label{remark4}
Under the assumption of Theorem \ref{borel-thm},
\begin{equation}
\label{Borel-gevrey-index}
r = -\frac{|b_i|}{w \cdot b_i}
\end{equation}
holds for any $i \not\in \sigma$,
where $\{b_i :\; i\notin \sigma \}$ is the basis of ${\rm{Ker}}\, A$
given in Remark \ref{remark3}.
Therefore if $u\in {\rm{Ker}}\, A \cap \mathbb{Z}^n$,
then $r w\cdot u$ is an integer.
\end{remark}

\begin{remark} The condition that ${A_w}$ contains $(1,1, \ldots, 1)$ is
assumed so that the Borel transformed series satisfies a regular
holonomic system {\rm \cite{hotta}}. Hence, the growth condition (ii) of
the Borel summability is satisfied because solutions of regular
holonomic systems satisfy a polynomial growth condition. Without this assumption, things become more complicated. See also
Section \ref{section:borel_transform_revisited}.\end{remark}

{\it Proof}\/.{\em (of Theorem \ref{borel-thm})}.
First of all, the open set $U$ in the theorem can be chosen as follows. There exist constants
$c_{ij}, c_i, p_i,m$ such that the series $\varphi_B(x,\zeta)$, which will be
defined in the proof below, converges when  $(x_1, \ldots, x_n,\zeta)$ belongs to the non empty open set $W=W' \cap(\cap_{j=1}^n (x_j\not= 0))$ where $W'$ is defined by the inequalities $\sum_{j} c_{ij} \log |x_j| + c_i \log |\zeta| <
p_i$, for $i=1, \ldots, m$. Such constants exist because $\varphi_B(x,\zeta)$
is a hypergeometric series which satisfies a regular holonomic
$A$-hypergeometric system ({\rm \cite{gkz1989}, \cite[Section 2.5]{SST}});
see also forthcoming Lemma \ref{borel-lemma2}. Since only
non-negative powers of $\zeta$ modulo an exponent appear in
$\varphi_B$, we may assume that $c_i > 0$ for $i=1,\ldots, m$. We may choose a non empty domain $U \subset \CC^n$ with  compact closure such that $U\times \{\zeta \in \CC | \, 0 < |\zeta |<\epsilon \} \subset W$ for some $\epsilon >0$.

To study the analytic properties of $\hat{\mathcal{B}}_{1/r}[\psi]$,
it is convenient to use
\begin{equation}
\label{solution-with-z}
\varphi(x, z)
:= \psi(x, t)\Big|_{t = z^r}
= \sum_{\ell = 0}^\infty C_{\ell} (x) z^{r (\ell + \gamma)}.
\end{equation}
Since $t^{-\gamma}\psi(x, t)$ is a formal power series in $t^{ w \cdot b_i}$,
$z^{-r\gamma} \varphi(x, z)$ does not contain any fractional powers in $z$
(cf. Remark \ref{remark4}).
Then  we have
\begin{equation}
\label{Borel-tr-of-formal-solution2}
\hat{\mathcal{B}}_1[\varphi](x, \zeta)
= \sum_{\ell = 0}^\infty
\frac{C_{\ell} (x)}{\Gamma(1+ r (\ell + \gamma))} \zeta^{r (\ell + \gamma)}
= \hat{\mathcal{B}}_{1/r}[\psi](x, \zeta^{r}).
\end{equation}
In what follows we simply write
$\varphi_B(x, \zeta)$
(resp., $\psi_B(x, \tau)$) instead of
$\hat{\mathcal{B}}_1[\varphi](x, \zeta)$
(resp., $\hat{\mathcal{B}}_{1/r}[\psi](x, \tau)$).

\begin{lemma}
\label{borel-lemma1} Assume condition \eqref{borel-expnent} holds for $\kappa=1/r$. For the power series $\varphi$ given in \eqref{solution-with-z}, we have
$$
\theta_{\zeta} \varphi_B
=
\hat{\mathcal{B}}_1[\theta_z \varphi] \,\,\, \mbox{ and } \,\,\,
\frac{\partial \varphi_B}{\partial\zeta}
=
\hat{\mathcal{B}}_1[z^{-1} \varphi].
$$
\end{lemma}

{\it Proof}\/.
The first relation follows from
\begin{align*}
\zeta\frac{\partial\varphi_B}{\partial\zeta}(x, \zeta)
&=
\sum_{\ell = 0}^\infty
\frac{C_{\ell} (x)}{\Gamma(1+ r (\ell + \gamma))} r(\ell + \gamma) \zeta^{r (\ell + \gamma)}
\\
&=
\hat{\mathcal{B}}_1
\left[
\sum_{\ell = 0}^\infty
C_{\ell} (x) r(\ell + \gamma) z^{r (\ell + \gamma)}
\right]
=
\hat{\mathcal{B}}_1
\left[
z \frac{\partial \varphi}{\partial z}
\right].
\end{align*}
We also have
\begin{align*}
\frac{\partial\varphi_B}{\partial\zeta}(x, \zeta)
&=
\sum_{\ell = 0}^\infty
\frac{C_{\ell} (x)}{\Gamma(1+ r (\ell + \gamma))} r(\ell + \gamma) \zeta^{r (\ell + \gamma) - 1}
\\
&=
\sum_{\ell = 0}^\infty
\frac{C_{\ell} (x)}{\Gamma(r (\ell + \gamma))}  \zeta^{r (\ell + \gamma) - 1}
\\
&=
\hat{\mathcal{B}}_1
\left[
\sum_{\ell = 0}^\infty
C_{\ell} (x) z^{r (\ell + \gamma)-1}
\right]
=
\hat{\mathcal{B}}_1
\left[
z^{-1} \varphi
\right].
\qquad
\QED
\end{align*}

\begin{lemma}
\label{borel-lemma2}
The formal power series $\varphi_B(x,\zeta)$ given in \eqref{Borel-tr-of-formal-solution2}
formally satisfies the hypergeometric system $H_{A_B}(\beta_B)$,
where
$$
A_B =
\begin{pmatrix}
A & 0\\
w & -1/r
\end{pmatrix},
\quad
\beta_B = \begin{pmatrix}\beta \\ 0 \end{pmatrix}.
$$
\end{lemma}

When the matrix
$A_B$
contains a rational entry, we regard the $\ZZ$-module generated
by the column vectors as the lattice to define the
$A$-hypergeometric system. For example, when
$A_B=\left(\begin{array}{ccc} 1&3&0\\1& 1&-1/2
\end{array}\right),$
$\beta_B=(\beta,0)$, the lattice is $\ZZ \times
\ZZ/2$ and the hypergeometric system is nothing but that for
$A_B=\left(\begin{array}{ccc} 1&3&0\\2& 2&-1
\end{array}\right)$ and $\beta_B=(\beta,0)$ for the lattice $\ZZ^2$.

{\it Proof}\/. It follows from Lemma \ref{borel-lemma1} and
relations $\theta_j \varphi = \theta_j \psi|_{t = z^r}$, $\theta_z
\varphi(x, z) = r (\theta_t \psi)|_{t = z^r}$ that
\begin{align*}
\left(\sum_{j = 1}^n a_{ij} \theta_j - \beta\right)\varphi_B
&=
\hat{\mathcal{B}}_{1}
\left[
\left(
\sum_{j = 1}^n a_{ij} \theta_j - \beta \right)\varphi
\right]
\\
&=
\hat{\mathcal{B}}_{1}
\left[
\left.
\left(
\sum_{j = 1}^n a_{ij} \theta_j - \beta \right)\psi
\right|_{t = z^r}
\right]
=0
\end{align*}
and
\begin{align*}
\left(\sum_{i = 1}^n w_i \theta_i - \frac{1}{r}\theta_{\zeta}\right)\varphi_B
&=
\hat{\mathcal{B}}_{1}
\left[
\left(\sum_{i = 1}^n w_i \theta_i - \frac{1}{r}\theta_z \right)\varphi
\right]
\\
&=
\hat{\mathcal{B}}_{1}
\left[
\left.
\left(\sum_{i = 1}^n w_i \theta_i - \theta_t \right)\psi
\right|_{t = z^r}
\right] = 0.
\end{align*}

Now we take vectors
$u = (u_1, \ldots, u_{n + 1})^T, v = (v_1, \ldots, v_{n + 1})^T  \in \mathbb{N}^{n+1}$
satisfying $A_B u =A_B v$.
By its definition, we obtain
$$
\dfrac{1}{r}(u_{n + 1} -v_{n + 1}) =
\sum_{i = 1}^n w_i (u_i - v_i) \in \mathbb{Z}.
$$
Without loss of generality we may assume one of $u_{n +1}$ and $v_{n + 1}$ is zero.
Under this assumption,
$u_{n + 1}/r$ and $v_{n + 1}/r$ are non-negative integers.
Furthermore
\begin{align*}
 &w_1 u_1 + \cdots + w_n u_n - \frac{1}{r} u_{n + 1}
= w_1 v_1 + \cdots + w_n v_n - \frac{1}{r} v_{n + 1}
\\
&\qquad\Longleftrightarrow\quad
w_1 u_1 + \cdots + w_n u_n + \frac{1}{r} v_{n + 1}
= w_1 v_1 + \cdots + w_n v_n + \frac{1}{r} u_{n + 1}
\end{align*}
holds. Therefore
$$
\begin{pmatrix} u_1 \\ \vdots \\ u_n \\ v_{n + 1}/r \end{pmatrix},
\begin{pmatrix} v_1 \\ \vdots \\ v_n \\ u_{n + 1}/r \end{pmatrix}
\in \mathbb{N}^{n + 1} \quad\text{and}\quad \tilde{A}
\begin{pmatrix} u_1 \\ \vdots \\ u_n \\ v_{n + 1}/r \end{pmatrix}
=
\tilde{A}
\begin{pmatrix} v_1 \\ \vdots \\ v_n \\ u_{n + 1}/r \end{pmatrix}.
$$
Hence
\begin{align*}
&
\big(\partial_1^{u_1} \cdots \partial_n^{u_n} \partial_{\zeta}^{u_{n + 1}}
- \partial_1^{v_1} \cdots \partial_n^{v_n} \partial_{\zeta}^{v_{n + 1}}\big)
 \varphi_B
\\
&\qquad
=
\hat{\mathcal{B}}_{1}
\left[
\big(\partial_1^{u_1} \cdots \partial_n^{u_n} z^{-u_{n + 1}}
- \partial_1^{v_1} \cdots \partial_n^{v_n} z^{-v_{n + 1}}\big)  \varphi
\right]
\\
&\qquad
=
\hat{\mathcal{B}}_{1}
\left[
z^{-(u_{n + 1} + v_{n + 1})}
\big(\partial_1^{u_1} \cdots \partial_n^{u_n} z^{v_{n + 1}}
- \partial_1^{v_1} \cdots \partial_n^{v_n} z^{u_{n + 1}}\big)
 \varphi
\right]
\\
&\qquad
=
\hat{\mathcal{B}}_{1}
\left[
z^{-(u_{n + 1} + v_{n + 1})}
\big(\partial_1^{u_1} \cdots \partial_n^{u_n} t^{v_{n + 1}/r}
- \partial_1^{v_1} \cdots \partial_n^{v_n} t^{u_{n + 1}/r}\big)
 \psi\,
\Big|_{t = z^r}
\right]
\\
&\qquad
= 0.
\end{align*}
Here we have used the second relation of Lemma \ref{borel-lemma1}.
This completes the proof. \QED

Under the assumption of Theorem \ref{borel-thm}, $H_{A_B}(\beta_B)$
is regular holonomic \cite{hotta}. Therefore, up to an non zero constant,  $\varphi_B(x, \zeta)$ is nothing but
a GKZ series solution of $H_{A_B}(\beta_B)$, and $\varphi_B(x,
\zeta)$ converges near $\zeta = 0$ if $x \in U$.
It also follows that the restriction of
$\varphi_B(x, \zeta)$ to $\{x = x^{0}\}$, which is a function of
$\zeta$, satisfies some
linear ordinary differential equation
$E(x^{0})$ of Fuchsian type.
Let ${\rm{Sing}}\,(x^{0})$ be the set  of
singular points of $E(x^0)$ except the origin and infinity and define
$\Theta(x^0) = \{\arg u;\, u \in {\rm Sing}\,(x^0)\}$.
For any $\theta \in \mathbb{R}$ with
$\theta\not\in\Theta(x^0)$ we have,
\begin{itemize}
\item[(i)]
$\varphi_B(x^0, \zeta)$ can be analytically continued
to a sector $S(\theta, \delta)$ with some small $\delta > 0$
since there is no singular points on $\{\zeta;\, \arg \zeta = \theta\}$.
\item[(ii)]
The Borel transform $\varphi_B(x^0,\zeta)$ has polynomial growth with respect to $\zeta$
in $S(\theta, \delta)$ since a singular point of $E(x^0)$
is a regular singular point.
\end{itemize}

Hence we conclude that $\varphi(x^0,z)$ is Borel summable
(i.e., 1-summable) in the direction $\theta$.
Since $\Theta(x^0)$ is a finite set for each fixed $x^0 \in U$,
$\varphi(x^0,z)$ is Borel summable in all but finitely many directions
for each fixed $x^0\in U$. We can consider a non empty open subset $U'\subset U$
such that the closure of $U'$ is compact and included in $U$. Then the complement in
$\RR$ of the union $\cup_{x\in U'} \Theta(x)$ contains an open interval. So for any
$\theta$ in this complement we can define ${\mathcal S}[\psi](x,t)$ using $[0, e^{i \theta})$
as a path of integration for any $x\in U'$. We may also avoid if necessary  the discriminant
of the leading term of the linear ordinary differential operator $E(x)$. We still write $U'=U$.

Because of the relation \eqref{Borel-tr-of-formal-solution2},
$\psi$ is $1/r$-summable if and only if  $\varphi$ is $1$-summable.
Therefore $\psi(x,t)$ is $1/r$-summable in all the directions in an open interval for all points $x\in U$.

Finally we show that the Borel sum of $\psi(x,t)$ is a solution
of the modified hypergeometric system $H_{A, w}(\beta)$
to finish the proof of Theorem \ref{borel-thm}.
Because of the relation
\begin{align*}
\mathcal{S}[\psi]\Big|_{t = z^r}
= \mathcal{L}^{\theta}_{1/r}\hat{\mathcal{B}}_{1/r}[\psi]\Big|_{t = z^r}
=\mathcal{L}^{\theta/r}_{1}\hat{\mathcal{B}}_{1}[\varphi]
= \mathcal{S}[\varphi],
\end{align*}
Lemma \ref{borel-lemma1} and Lemma \ref{borel-lemma2},
it is enough to prove
(we omit the direction $\theta/r$ in the following.)

\begin{lemma}
$$
z \frac{\partial}{\partial z} \mathcal{L}_1[\varphi_B]
= \mathcal{L}_1 \left[\zeta \frac{\partial\varphi_B}{\partial \zeta}\right],
\quad
z^{-1} \mathcal{L}_1[\varphi_B] =
\mathcal{L}_1\left[\frac{\partial\varphi_B}{\partial \zeta}\right].
$$
\end{lemma}

{\it Proof}\/.
To begin with, we give a proof in the case when $\Re (r\gamma) > 0$ or $\gamma=0$.
By differentiating under the integral sign we obtain
\begin{align*}
z \frac{\partial}{\partial z} \mathcal{L}_1[\varphi_B]
&=
\int^{\infty}_0 \dfrac{\zeta}{z} e^{-\zeta/z} \varphi_B(x, \zeta) \dfrac{d\zeta}{z}
-
\int^{\infty}_0  e^{-\zeta/z} \varphi_B(x, \zeta) \dfrac{d\zeta}{z}
\\
&=
-
\int^{\infty}_0
\frac{\partial}{\partial \zeta}
(\zeta e^{-\zeta/z}) \cdot  \varphi_B(x, \zeta) \dfrac{d\zeta}{z}.
\end{align*}
By integral by parts, this equals
\begin{align*}
-
\left[
\zeta e^{-\zeta/z} \cdot  \dfrac{\varphi_B(x, \zeta) }{z}
\right]_{\zeta = 0}^{\infty}
+
\int^{\infty}_0
e^{-\zeta/z} \cdot \zeta \frac{\partial \varphi_B}{\partial \zeta}(x, \zeta) \dfrac{d\zeta}{z}.
\end{align*}
Since the boundary terms vanish
(the boundary term coming from infinity vanishes because
of the growth estimate of $\varphi_B$), the first relation follows.
In a similar manner, we obtain
\begin{align*}
\mathcal{L}_1\left[\frac{\partial\varphi_B}{\partial \zeta}\right]
&=
\int^{\infty}_0 e^{-\zeta/z} \frac{\partial\varphi_B}{\partial\zeta}(x, \zeta) \dfrac{d\zeta}{z}
\\
&=
\left[
e^{-\zeta/z}
\frac{\varphi_B(x, \zeta)}{z}
\right]_{\zeta = 0}^{\infty}
+ \frac{1}{z}
\int^{\infty}_0 e^{-\zeta/z} \varphi_B(x, \zeta) \dfrac{d\zeta}{z}.
\end{align*}
The growth estimate of $\varphi_B$ at infinity
and the behavior of $\varphi_B$ near the origin guarantee that
the boundary terms vanish.
This proves the second relation.

When $\Re (r\gamma) \leq  0$ and $\gamma$ satisfies \eqref{borel-expnent},
we use \eqref{def-Borel-sum2} as the definition of the Borel sum,
and the same argument works.
In this case all of the boundary terms come from infinity
and  they vanish.
\QED This finish the proof of Theorem \ref{borel-thm}. \QED

In the preceding proof, we have shown that $\psi_B(x,\tau)$ is of polynomial growth in $\tau$ for $x\in U$. Therefore,
for an arbitrarily small positive $c_2$
we can find $c_1 > 0$ for which \eqref{borel-summability-ii} holds
with $f = \psi$, $\kappa = 1/r$ for all $x\in U$ after eventually replacing $U$ by a smaller open set. Therefore, the condition \eqref{def-Borel-sum-domain} guarantees
the following Corollary.
\begin{corollary}\label{corollary-asymptotic}
The Laplace integral \eqref{def-Borel-sum} of the Borel transform
$\psi_B(x,\tau)$ converges  in $S(\theta, r \pi)$ if $\theta
\not\in\Theta(x)$. In particular, we can set $t = 1$ in the
expression ${\mathcal S}[\psi](x,t)$ of the Borel sum of $\psi(x,t)$
if $\Theta(x)$ does not contain $0$.
The series ${\mathcal S}[\psi](x,1)$  gives a holomorphic solution of $H_A(\beta)$ and the
formal series $\psi(x, 1)$ also expresses the asymptotic expansion of
the solution ${\mathcal S}[\psi](x,1))$ along a curve
$x(t)= (t^{w_1} x^0_1, \ldots, t^{w_n}
x^0_n)$ as $t\rightarrow 0$ for $x^0\in U$. That is $$|\mathcal{S}[\psi](x(t) ,1)- \psi_N (x(t),1)| \leq C K^N |t|^N\Gamma (1+ r N)$$ for all $N\geq 0$ for some constants $C>0,K>0$,  where $\psi_N(x,t):=\sum_{\ell < N} C_\ell(x)t^\ell$.
\end{corollary}

\begin{remark} A comment on the statement of previous Corollary.
A priori the constants $c_1, c_2$ of the estimate \eqref{borel-summability-ii} for $\psi(x,t)$ and $\kappa=1/r$ depend on $x\in U$. Since $\varphi_B(x^0,\zeta)$ is holomorphic at $(x^0,0)$ for $x^0\in U$, we can find uniform constants $c_1, c_2$ for $x$ in a small neighborhood of each $x^0$ in $U$ with $\theta \not\in \Theta(x^0)$. Notice also that, from Lemma \ref{formal-solution-w-perturbation},  $\psi(x,t)=\phi_v(t^wx)=\psi(t^wx,1)$ and that  $\psi(x,1)=\phi_v(x)$ is a (possibly Gevrey) solution of $H_A(\beta)$.
We also have the following formal equalities: $\psi_B(x,\tau)=\psi_B(\tau^w x,1)$ and ${\mathcal S}[\psi](x,t)={\mathcal S}[\psi](t^wx,1)$.

Then, fixing $x^0\in U$, the estimate \eqref{eq:asymptotic_expansion} evaluated at points  $(t^w x^0,1)$ proves that the formal expression $\psi(t^w x^0,1)$ is an asymptotic expansion of $t^{-\gamma} \mathcal{S}[\psi](t^w x^0,1)$ when $t \rightarrow 0$.
\end{remark}

\begin{remark}
\label{remark6}
Since the singular points ${\rm{Sing}}\,(x)$ of $E(x)$ depend
on $x$, one of them may meets the path of integration
of the Borel sum if $x$ moves.
In that case we obtain the analytic continuation of the Borel sum
by deforming the path of integration.
This is closely related to the Stokes phenomenon.
\end{remark}

\begin{example} \rm  \label{example12}
Put $A=(1,2)$. Then, the image of ${\bar A}^T$ is $\RR^2$.
Let us take $w=(0,1)$.
A formal solution of the modified system $H_{A,w}(\beta)$ is
$$ \psi(x,t) = x_1^\beta \sum_{m=0}^\infty
  \frac{[\beta]_{2m}}{m!} \left( \frac{x_2}{x_1^2} \right)^m t^m.
$$
If $\beta \not \in \NN$, the Gevrey index $s= r+1$ of $\psi(x,t)$
along $t=0$ is $s=2$. We set $z = t$ and $\varphi(x, z) = \psi(x, z)$.
The Borel transform $\varphi_B$ is
\begin{eqnarray*}
x_1^\beta \sum_{m=0}^\infty
  \frac{[\beta]_{2m}}{m!} \left( \frac{x_2}{x_1^2} \right)^m \frac{\zeta^m}{m!}
 = x_1^\beta \cdot {}_2F_1(-\beta/2,(-\beta+1)/2,1;4 x_2 \zeta/x_1^2).
\end{eqnarray*}
The domain $U$ may be defined by $\{ (x_1, x_2) \,|\, -2 \log |x_1|
+ \log |x_2| < -1, |x_1| < 1 \}$. The series $\varphi_B$ satisfies
the $A$-hypergeometric system with $A_B= \left( \begin{array}{ccc}
 1 & 2 & 0 \\
 0 & 1 & -1 \\
\end{array} \right)
$
and $\beta_B=(\beta,0)^T$.
The equation $E(x) = E(x_1, x_2)$ in the proof is
$$
\left[
(4 x_2 \zeta^2 - {x_1}^2\zeta)
\left(\frac{\partial}{\partial\zeta}\right)^2
+
\big((-4\beta + 6) x_2 \zeta - {x_1}^2 \big)\frac{\partial}{\partial\zeta}
+ (\beta^2 - \beta)x_2
\right] \varphi_B(x, \zeta) = 0,
$$
and
$$
{\rm{Sing}}\,(x) = \left\{{x_1}^2/(4x_2)\right\},
\quad
\Theta(x) = \{2\arg x_1 - \arg x_2\}.
$$
The Borel sum of $\psi(x, t)$ is
$$
\dfrac{x_1}{t}
\int^{e^{i\theta}\cdot \infty}_0 e^{-\tau/t}
{}_2F_1(-\beta/2,(-\beta+1)/2,1;4 x_2 \tau/x_1^2)
d\tau
$$
and, for each $x \in U$,
$\varphi$ (and hence $\psi$) is Borel summable in all directions except
the angle $\theta = 2\arg x_1 - \arg x_2$.

The series $\psi(x,1)$ can be regarded as an asymptotic expansion
of a solution
of the original $A$-hypergeometric system for $A=(1,2)$ and $\beta$
from (\ref{eq:asymptotic_expansion}) and Corollary \ref{corollary-asymptotic}.
In other words, we have
$$ x_1^\beta
\int_0^{e^{i\theta}\infty}
e^{-\tau}
{}_2 F_1\left( \frac{-\beta}{2},\frac{-\beta+1}{2}, 1;
 \frac{4x_2 \tau}{x_1^2} \right)
 d\tau
\sim \psi(x,1),$$
which is a well-known asymptotic expansion.
\end{example}

\begin{example} \rm  \label{example1356}
Put $A=(1,3,5,6) $.
Then, the image of ${\bar A}^T$ is
$$ \{ (w_1,w_2,w_3,w_4) \,|\, w_2-3 w_3+2w_4 = 0,
  w_1 - 5 w_3 + 4 w_4 = 0 \}.
$$
Let us take $w=(-4,-2,0,1)$ in the ${\rm Im}\, {\bar A}^T$.
Then, the initial ideal ${\rm in}_{(-w,w)}(H_A(\beta))$
is generated by
$\{  \partial_{{x}_{2}}, \partial_{{x}_{3}}, \partial_{{x}_{4}} \}$
and
$ \theta_1 + 3 \theta_2 + 5 \theta_3 + 6 \theta_4-\beta$.
The rank of this system is $1$
and the solution of this system is spanned by
$x_1^\beta$.
We extend this solution to a series solution of $H_A(\beta)$.
The solution can be written as
\begin{equation}
\phi_v = \sum_{u \in N_v} \frac{[v]_{u_-}}{[v+u]_{u_+}} x^{v+u}=
x_1^{\beta} \big(
  1
+ \frac{\beta(\beta-1)(\beta-2)}{1!} x_1^{-3} x_2 + \cdots \big)
\label{eq:ex1356series2}
\end{equation}
where $v = (\beta,0,0,0)$.
The corresponding series solution of the modified system is
\begin{eqnarray}\label{eq:ex1356seriesm}
&& \psi_v (x ,t) = \phi_v (t^{-4} x_1 , t^{-2} x_2
, x_3 , t x_4 ) =x_1^{\beta} t^{-4\beta}\big(
  1
 + \frac{\beta(\beta-1)(\beta-2)}{1!} x_1^{-3} x_2 t^{10}
+ \cdots
\big)
\end{eqnarray}
The Gevrey index $r+1=1/\kappa+1$ of the series $\psi_v(x,t)$ with
respect to $t=0$ is $1/5+1$ if $\beta$ is very generic, see Remark
\ref{remark4}. Apply the Borel transformation (\ref{Borel-tr-of-formal-solution2}).
The transformed series
\begin{eqnarray}\label{eq:ex1356seriesb}
&&x_1^{\beta} \zeta^{-4\beta/5}\big(
  \frac{1}{\Gamma(1-\frac{4}{5}\beta)}
+ \frac{\beta(\beta-1)(\beta-2)}{1!} x_1^{-3} x_2 \frac{\zeta^{2}}{\Gamma(1-\frac{4}{5}\beta+2)}  \\
&+& \nonumber
  \left(
  \frac{\beta(\beta-1)\cdots (\beta-5)}{2!} x_1^{-6} x_2^2
+ \frac{\beta(\beta-1)\cdots (\beta-4)}{1!} x_1^{-5} x_3^1
  \right) \frac{\zeta^{4}}{\Gamma(1-\frac{4}{5}\beta+4)}
+ \cdots
\big)
\end{eqnarray} satisfies the $A$-hypergeometric system
associated to the matrix
\begin{equation}
A_B =
\left(
\begin{array}{ccccc}
 1 & 3 & 5 & 6 & 0\\
-4 &-2 & 0 & 1 & -5 \\
\end{array}
\right)
\end{equation}
and $\beta_B = (\beta,0)$.
It follows from the condition on the weight vector $w$ that this
$H_{A_B}(\beta_B)$ is a regular holonomic system.
The series (\ref{eq:ex1356seriesb}) can be obtained
by taking the $(-u,u)$-initial ideal of $H_{A_B}(\beta_B)$
with respect to the weight vector $(w,0)$ and $(0,0,1,2,0)$
as the tie breaking weight vector \cite[Chapters 2,3]{SST}.

The kernel element $\ell$ of $A_B$ as a map from $\RR^5$ to $\RR^2$
is parametrized as
$$ \ell_1 = \ell_3 + \frac{3}{2} \ell_4 - \frac{3}{2} \ell_5,
   \ell_2 = -2 \ell_3 - \frac{5}{2} \ell_4 + \frac{1}{2} \ell_5.
$$
Since we sum the series on
$ \ell_i \geq 0 $, $i=2,3,4,5$ and $\ell \in \ZZ^5$,
these lattice elements $\ell$ can be parametrized as
$$ \ell_5 = 2m+e_p, \ell_4 = 2n+e_p, \ell_3 = k,
   \ell_2 = -2k-5n+m-2e_p, \ell_1 = k+3n - 3m
$$
$$ m,n \in {\NN}  \ \mbox{and} \ 2k + 5n + 2e_p \leq m $$
where $e_p = 0$ or $1$.
Let us introduce the following hypergeometric series
$$ \sum_{m,n,k \geq 0, 2k+5n+2e_p \leq m}
  \frac{ z_3^k z_4^{2n+e_p} z_5^{2m+e_p}}
 { (b)_{k+3n-3m} (-2k-5n+m-2e_p)! k! (2n+e_p)! (a)_{2m+e_p}}
$$
depending on an integer $e_p$.
We denote it by $F_{\rm odd}(a,b; z_3,z_4,z_5)$ when $e_p=1$
and by $F_{\rm even}(a,b; z_3,z_4,z_5)$ when $e_p=0$.
Then, the series \eqref{eq:ex1356seriesb} is expressed as
$$ \frac{x_1^\beta \zeta^{-4\beta/5}}{\Gamma(1-\frac{4}{5} \beta)}
   F(x_1 x_2^{-2}x_3,x_1^{3/2} x_2^{-5/2}x_4,x_1^{-3/2}x_2^{1/2}\zeta) $$
where
$$F(z) =  F_{\rm odd}(1-\frac{4}{5}\beta,\beta+1;z)
         +F_{\rm even}(1-\frac{4}{5}\beta,\beta+1;z).
$$
It follows from \eqref{eq:asymptotic_expansion} that
the series \eqref{eq:ex1356seriesm}  is an asymptotic expansion of
$$\int_0^{e^{i\theta} \infty}
  e^{-(\tau/t)^5}
 \frac{x_1^\beta \tau^{-4\beta/5}}{\Gamma(1-\frac{4}{5}\beta)}
 F(x_1 x_2^{-2}x_3,x_1^{3/2} x_2^{-5/2}x_4,x_1^{-3/2}x_2^{1/2}\tau)
  \frac{5 \tau^4}{t^5} d\tau .
$$
\end{example}

\section{Borel transformation revisited}
\label{section:borel_transform_revisited}

In this section we assume that the $\QQ$-row span of the matrix $A$ does
not contain the vector $(1,\ldots, 1)$ and we see that by studying the irregularity of
the modified $A$-hypergeometric system along $T$ we can give
an analytic meaning to the Gevrey series solutions of the
$A$-hypergeometric system $\cM_A (\beta)$ along coordinate varieties
constructed in \cite{Fer}. More precisely, we prove  (Remark
\ref{w-in-borel-revisited} and Theorem \ref{theorem-asymptotic-expansion}) that these Gevrey solutions of $\cM_A(\beta)$ are asymptotic expansions of
certain holomorphic solutions of $\cM_A (\beta)$ when some conditions are satisfied.

Let us start with an observation about the assumption of Theorem
\ref{borel-thm} that the $\QQ$-row span of the matrix $A$ does not contain
the vector $(1, \ldots, 1)$, but the one of $A_w$ does. Since adding to $w$ a linear combination of the
rows of $A$ does not change the ideal $H_{A,w ,\alpha }(\beta)$ if
we accordingly change $\alpha$, we can assume without loss of generality that $w = \lambda (1, \ldots ,1)$
for some non zero $\lambda\in \ZZ$. On the other hand, for $\lambda > 0$ the modified system $\cM_{A ,w , \alpha}
(\beta )$ is regular along $T$ by Proposition
\ref{regularity-along-T} and Remark \ref{regularity-along-T}. Thus,
we may further assume that $w = (-\kappa, \ldots, -\kappa)$ for some
$\kappa\in\ZZ_{>0}$. For this $w$ each of the formal series
constructed in the proof of Theorem \ref{theorem-dimension-1} is a
Gevrey series along $T$ with order $s=r+1=1+ 1/\kappa$ multiplied by
a term $t^{\gamma}$, $\gamma \in \CC$. This sort of formal series
solutions along $T$ of the modified system $\cM_{A ,w , \alpha}
(\beta )$ includes series of the form $t^{-\alpha} \phi ( t^{w_1 }
x_1 ,\ldots , t^{w_n}x_n )$ where $\phi (x)$ is a Gevrey solution of
$\cM_A (\beta )$ along a coordinate variety $Y$ of low dimension which is
not Gevrey along a coordinate variety of greater dimension (see the proof of Proposition
\ref{Gevrey-for-special-weight-vector}).

Recall that a vector $v$ is said to be associated with a simplex $\sigma$ if $v= v^{\mathbf{k}}$
is such that $v_i = k_i \in \NN $ for all $i \notin \sigma$ and $A v
= \beta$. Let us also recall that $\Delta_A$ is the convex hull of
the columns of $A$ and the origin while $\operatorname{conv}(A) $ is
the convex hull of the columns of $A$. We have the following

\begin{proposition}\label{Gevrey-for-special-weight-vector}
Take $w = (-\kappa , \ldots , -\kappa )$ with $\kappa\in\ZZ_{>0}$ and assume
that $\beta$ is very generic. Let $\sigma$ be a $(d-1)$-simplex of $A$
and $v$ a vector associated with $\sigma$. Then, up to
multiplication by a term $t^{\gamma}$, for some $\gamma \in \CC$, the series
$t^{-\alpha} \phi_{v} ( t^{w_1 } x_1 ,\ldots , t^{w_n} x_n )$ is a
formal power series along $T$ (and in such a case, it is a Gevrey solution of
$\cM_{A, w ,\alpha }(\beta)$ with index $s'=r'+1=1+1/\kappa$) if and
only if $\sigma$ is contained in a facet of $\operatorname{conv}(A)$
that is not a facet of $\Delta_A$.
\end{proposition}

{\it Proof}\/. From \cite[Theorem 3.11]{Fer}, $\phi_v$ is a Gevrey solution of
$\cM_A (\beta)$ along $Y=\{x_i = 0 :\: |A_{\sigma}^{-1}a_i|> 1 \}$
with order $s=r+1= \max_{i} \{|A_{\sigma}^{-1}a_i| \}$. In fact, for
any simplex $\sigma$ of $A$ we can construct
$\operatorname{vol}(\sigma)$ linearly independent Gevrey solutions
as before (see \cite[Remark 3.6]{Fer}).

When
$\beta$ is very generic these series have Gevrey index $s=r+1=
\max_{i} \{|A_{\sigma}^{-1}a_i| \}$ and they are of the form
$\phi_v$ for $v$ associated with $\sigma$.

The series $t^{-\alpha} \phi_{v} ( t^{w_1 } x_1 ,\ldots , t^{w_n}
x_n )$ is annihilated by $H_{A,w,\alpha}(\beta)$ because $\phi_v$ is
annihilated by $H_A (\beta)$. Moreover, monomials appearing in
$\phi_v$ are of the form $x^{v+u}$ for $u$ integer vectors in an
affine translate of the positive span of the columns of the matrix
$B_{\sigma}$. Thus $t^{-\alpha} \phi_{v} ( t^{w_1 } x_1 ,\ldots ,
t^{w_n} x_n )$ is a formal series along $T$ (up to multiplication by
$t^{\alpha -w v}$) if and only if all the exponents of $t$ belong to
$\alpha -w v + \NN$. This happens if and only if $w u \in \NN$ for
all the columns $u$ of $B_{\sigma}$. For $w = (-\kappa, \ldots,
-\kappa)$ with $\kappa\in\ZZ_{>0}$ the scalar product of $w$ and a
column of $B_{\sigma}$ is given by $\kappa (|A_{\sigma}^{-1}a_i|-1)$
for $i\notin \sigma$ and this product is nonnegative for all
$i\notin \sigma$ if and only if $|A_{\sigma}^{-1}a_i|\geq 1$ for all
$i\notin \sigma$. This is equivalent to say that $\sigma$ is a
simplex of $A$ such that all the columns of $A$ are either in the
hyperplane $H_{\sigma}$ passing through the columns of $A$ indexed
by $\sigma$ or in the corresponding single open half space not
containing the origin. Equivalently, $\sigma$ is contained in a
facet of $\operatorname{conv}(A) $ that is not a facet of
$\Delta_A$. \QED

The Gevrey series solutions constructed in \cite{Fer} are of the form $\phi_v$ for $v$ associated with a
suitable simplex $\sigma$ of $A$. In order to interpret some of these Gevrey series as an asymptotic expansion
of a solution of $M_A(\beta)$ via the modified $A$-hypergeometric system, we consider a vector $w\in \ZZ^n$
with the following coordinates:
\begin{equation}
 w_i =\left\{\begin{array}{ll}
 |\operatorname{det}(A_{\sigma})|(|A_{\sigma}^{-1}a_i|-1) & \mbox{ if } |A_{\sigma}^{-1}a_i|> 1  \\
0  & \mbox{ otherwise }
 \end{array}\right. \label{sigma-weight-vector}
\end{equation} up to addition with a linear combination of the rows of $A$.
Notice that when $\sigma$ is contained in a facet of
$\operatorname{conv}(A)$ that is not a facet of $\Delta_A$, this
vector $w$ verifies the assumptions of Theorem \ref{borel-thm}. More
precisely, $|\operatorname{det}(A_{\sigma})|(1 ,\ldots ,1)+w$ is a
linear combination of the rows of $A$.

\begin{remark}\label{two-hyperplanes}
Notice that for $w$ given by  (\ref{sigma-weight-vector}) all the
columns $\wa_i$ of the matrix $\wA$ except for $\wa_{n+1}$ are contained in
the union of at most two hyperplanes. If we take
coordinates $(y,y_{d+1})$ in $\RR^{d} \times \RR$, then the
hyperplane $\{ y_{d+1} = 0 \}$ contains all the columns $\wa_i$ such
that $|A_{\sigma}^{-1}a_i|\leq 1$ and the hyperplane $\{
|A_{\sigma}^{-1}y|-\frac{1}{|\operatorname{det}(A_{\sigma})|}
y_{d+1} = 1 \}$ contains
$-|\operatorname{det}(A_{\sigma})|\wa_{n+1}$ and all the columns
$\wa_i$ such that $|A_{\sigma}^{-1}a_i|\geq 1$. In particular, the
intersection of these two hyperplanes contains all the columns
$\wa_i$ such that $|A_{\sigma}^{-1}a_i|=1$ (for example, all the
columns $\wa_i$ for $i\in \sigma$). We also notice that the
points $\{\wa_i : \: i=1,\ldots ,n\}\cup
\{-|\operatorname{det}(A_{\sigma})|\wa_{n+1}\}$ belong to the same
hyperplane if and only if $\sigma$ is contained in a facet of
$\operatorname{conv}(A)$ that is not a facet of $\Delta_A$. Notice
that if $\sigma$ is contained in a facet of $\Delta_A$ and $v$ is
associated with $\sigma$, then the corresponding $w$  given by (\ref{sigma-weight-vector}) is $0$ and
$\phi_v(x)$ is convergent.
\end{remark}

In particular, Remark \ref{two-hyperplanes} proves the following:

\begin{lemma}
If $\sigma$ is a $(d-1)$-simplex of $A$ not contained in a facet of $\Delta_A$ and $w$ is given by
(\ref{sigma-weight-vector}) then $s'= r'+1=1 + 1/
|\operatorname{det}(A_{\sigma})|$ is a slope of
$\cM_{A,w,\alpha}(\beta)$ along $T$, for any $\alpha$.
\end{lemma}

\begin{proposition}\label{Gevrey-for-sigma-weight-vector}
Let $\sigma$ be a $(d-1)$-simplex of $A$ not contained in a facet of
$\Delta_A$ and consider $w$ given by (\ref{sigma-weight-vector}) and
$\beta$ very generic. For any vector $v$ associated with $\sigma$ we
have that, up to multiplication by $t^{\alpha -w v}$, the series
$\psi (x ,t)=t^{-\alpha} \phi_{v} ( t^{w_1} x_1 ,\ldots , t^{w_n}
x_n )$ is a Gevrey solution of $\cM_{A, w ,\alpha }(\beta)$ with
index $s'=r'+1=1+1/|\operatorname{det}(A_{\sigma})|$ along $T$
at any point of $T\cap U_{\sigma , R}$ for some $R>0$, where
$U_{\sigma ,R}=\{(x, t ) \in \CC^{n}\times \CC :\;  | x_j t^{w_j} |<
R |x_{\sigma}^{A_{\sigma}^{-1} a_j}| \mbox{ if } j \notin \sigma
\mbox{ and } |A_{\sigma}^{-1} a_j | \geq 1 \}\cap \{ x_i \neq 0 :\;
i \in \sigma \}$ .
\end{proposition}

{\it Proof}\/. By \cite[Theorem 3.11]{Fer} we have that
$\phi_v(x)$ is a Gevrey solution of $\cM_A (\beta)$ with index $s=r+1=
\max_{i} \{|A_{\sigma}^{-1}a_i| \}$) along $Y=\{x_i = 0 :\:
|A_{\sigma}^{-1}a_i|> 1 \}$ at any point of $Y\cap\{x \in \CC^{n}
:\;  | x_j  |< R |x_{\sigma}^{A_{\sigma}^{-1} a_j}| \mbox{ if } j
\notin \sigma \mbox{ and } |A_{\sigma}^{-1} a_j | = 1 \}\cap \{ x_i
\neq 0 :\; i \in \sigma \}$,
for some $R>0$.

It is clear from (\ref{sigma-weight-vector}) that $w\in \NN^n$ and $w_j =0$ for all
$j\in \sigma$. Hence, for any exponent $v+u$ in the series $\phi_v(x)$
the corresponding exponent of $t$ in the series $\psi (x ,t)$ is $-\alpha + w (v+u)=
-\alpha +\sum_{j\notin \sigma} w_j (v_j + u_j) \in -\alpha + \NN$ because $v_j \in \NN$ for
all $j\notin \sigma$ and $u\in N_v$. We conclude the proof by using Remark
\ref{remark3}. \QED

In analogy with Section \ref{section:Borel} we denote $\psi_B
(x,\tau )=\hat{\mathcal{B}}_{1/r'}[\psi](x, \tau)$, which defines a
holomorphic function at any point in $U_{\sigma , R}$ for some $R>0$
by Proposition \ref{Gevrey-for-sigma-weight-vector}. We also denote
$\varphi (x ,z)=\psi (x, z^{r'} )$ and hence $\varphi_B (x,\zeta
)=\hat{\mathcal{B}}_{1}[\varphi](x, \zeta )=\psi_B ( x ,
\zeta^{r'})$ is convergent at points in the open set
\begin{equation}
U_{\sigma ,R} ' =\{ ( x ,\zeta ) \in \CC^{n}\times \CC :\; \zeta = \tau^{|\operatorname{det}(A_{\sigma}) |},\; (x ,\tau )\in U_{\sigma , R} \} \label{open-set}
\end{equation}

Moreover, the series $\varphi_B (x, \zeta)$ is a holomorphic
solution of $H_{A_B}(\beta_B)$ (in the variables $(x,\zeta)$), where
$$
A_B =
\begin{pmatrix}
A & 0\\
w & -\kappa
\end{pmatrix},
\quad \beta_B = \begin{pmatrix}\beta \\ \alpha \end{pmatrix} \quad {\mbox { with  }} \kappa =|\operatorname{det}(A_{\sigma})|.$$

\begin{remark}
For $w$ given by (\ref{sigma-weight-vector}) the hypergeometric
system $H_{A_B}(\beta_B)$, can have slopes along $\zeta =\infty$ for
all $\beta , \alpha \in \CC$ (see Example
\ref{example-slope-infinity-2}). However, see Proposition
\ref{bounded-growth} where we point out a property of the
solution $\varphi_B(x,\zeta)$.
\end{remark}

\begin{example}\label{example-slope-infinity-2}
Let us consider the matrix $$A=\left(\begin{array}{cccc}
                                2 & 0 & 1 & 3\\
                                0 & 1 & 1 & 2\end{array}\right),$$ $\beta\in \CC^2$, $\alpha \in \CC$, the simplex $\sigma =\{1, 3\}$ and $w=(0,0,0,3)$ given by
(\ref{sigma-weight-vector}). Notice that $\det (A_{\sigma})=2$ and
hence $$A_B=\left(\begin{array}{ccccc}
                                2 & 0 & 1 & 3 & 0\\
                                0 & 1 & 1 & 2 & 0\\
                                0 & 0 & 0 & 3 & -2\end{array}\right)$$ Using Corollary \ref{slopes-hypergeometric-at-infinity} we have
that $s=1+r= 1 + 1/3$ is a slope of $M_{A_B} (\beta_B )$ along $\zeta =\infty$.
\end{example}

\begin{lemma}\label{technical-lemma}
Assume $\Phi_{A_{\eta}}^{F,d-1}\subseteq \Phi_A^{F,d-1}$ where
$\eta=\{ i :\; |A_{\sigma }^{-1}a_i | \geq 1 \}$. Let
$\widetilde{\sigma}\subseteq \tau \in \Phi_{A_{\eta}}^{F, d-1}$ be a
simplex, $\widetilde{v}\in \CC^{n+1}$ a vector associated with $\widetilde{\sigma}\cup \{n+1\}$ (i.e. $A_B \widetilde{v}=\beta_B$ and $\widetilde{v}_i \in \NN$ for all $i\notin\widetilde{\sigma}\cup \{n+1\}$). The
series $\phi_{\widetilde{v}}(x ,\zeta)$ converges at points $(x ,\zeta )\in U \times \{ \zeta :\; |\zeta |>R '\}$
and for arbitrarily small $c_2>0$ we can choose $c_1>0$ such that $|\phi_{\widetilde{v}}(x ,\zeta)|\leq c_1 \operatorname{exp}(c_2 |\zeta |)$.
\end{lemma}

\begin{remark} The previous condition on the $(A,F)$-umbrella holds for any $d\times n$ matrix $A$ with $d=1$ or $n-1=d$.
\end{remark}

{\it Proof}\/. Notice that $\eta \cup \{n+1\}$ is a facet of the
$(A_B , F)$-umbrella and $\widetilde{\sigma}\cup\{n+1\}$ is a
simplex of $A_B$ contained in $\eta \cup \{n+1\}$. In particular we
know that $\phi_{\widetilde{v}}$ is convergent in certain open set.
Let us denote $\kappa =|\det (A_{\sigma })|$. If $\{ b_i :\; i\notin
\widetilde{\sigma}\cup \{n+1\}\}$ is the basis of $\ker (A_B )$
associated with $\widetilde{\sigma}\cup\{n+1\}$ then
the coordinate sum
of $b_i$ is
$$|b_i|=\left\{\begin{array}{cc}
               0 & \mbox{ if } i\in \eta \setminus \widetilde{\sigma}\\
               -|A_{\widetilde{\sigma}}^{-1}a_i| +1 - \frac{1}{\kappa} w_{\widetilde{\sigma}}
A_{\widetilde{\sigma}}^{-1}a_i & \mbox{ if } i\notin \eta
\end{array}\right.$$
Let us denote by $(b_i)_A$ the vector given by the first $n$ entries of $b_i$.
Since $\widetilde{\sigma}\cup \{n+1\} \subseteq \eta \cup \{n+1\}\in\Phi_{A_B}^{F,d-1}$,
we have that $|b_i |> 0$ for all $i\notin \eta $ and that the series $\phi_{\widetilde{v}}$
defines a multivalued holomorphic function in the open set
$\{(x,\zeta ):\; |x^{(b_i)_A} \zeta^{(b_i)_{n+1}} |< R, \; i\in\eta
\setminus \widetilde{\sigma}\}\cap \{ x_j \neq 0 :\;
j\in \widetilde{\sigma} \}$ for some $R>0$. The fact that
$\widetilde{\sigma}\subseteq\tau \in
\Phi_{A_{\eta}}^{F,d-1}\subseteq \Phi_A^{F,d-1}$ guarantees that
$-|A_{\widetilde{\sigma}}^{-1}a_i| +1 \geq 0$ for all $i\in \eta$
and $-|A_{\widetilde{\sigma}}^{-1}a_i| +1 >0$ for all $i\notin
\eta$. Thus, if $i\in \eta \setminus \widetilde{\sigma}$ the last
coordinate of $b_i$ is $(b_i)_{n+1}=|b_i|-|(b_i)_A|=
-1+|A_{\widetilde{\sigma}}^{-1}a_i|\leq 0$ while if $i\notin \eta$
the last coordinate of $b_i$ can be positive.

However, if $(b_i )_{n+1}>0$ for some $i\notin \eta$, we still have
that $|b_i|> (b_i)_{n+1}= -\frac{1}{\kappa} w_{\widetilde{\sigma}}
A_{\widetilde{\sigma}}^{-1}a_i$. In this case there exist $K_1 , K_2>0$ such that:
$$\sum_{m\geq 0} (m !)^{-|b_i|} |x^{(b_{i})_A}\zeta ^{(b_i)_{n+1}}|^m \leq K_1
\operatorname{exp}(K_2 |x^{(b_{i})_A}|^{1/|b_i|}|\zeta|^{(b_i)_{n+1}/|b_i|}) $$  where $(b_i)_{n+1}/|b_i|<1$.

On the other hand, if $\widetilde{v}_i = k_i$ for $i\notin \widetilde{\sigma}\cup \{n+1\}$,
it can be shown by using standard estimates on $\Gamma$-functions
(see e.g. \cite[Proposition 1, Section 1.1]{gkz1989}, \cite[Lemma 1]{ohara-takayama} and \cite[Lemma 3.8.]{Fer}),
that there exists $C_1, C_2>0$ such that

\begin{align*}
&|\phi_{\widetilde{v}}(x ,\zeta)|\leq C_1
|x^{A_{\widetilde{\sigma}}^{-1}\beta}\zeta^{(-\alpha
+w_{\widetilde{\sigma}} A_{\widetilde{\sigma}}^{-1}\beta)/\kappa}|
\sum_{k+m\in \NN^{n-d}} \frac{C_2^{\sum k_i + m_i } |x^{\sum_i (k_i
+ m_i ) (b_i)_A} \zeta^{\sum (k_i + m_i )
(b_i)_{n+1}}|}{\prod_{i\notin \widetilde{\sigma}\cup \{n+1\}} (k_i +
m_i )!^{|b_i|}}=
\\
&=
C_1 |x^{A_{\widetilde{\sigma}}^{-1}\beta}\zeta^{(-\alpha
+w_{\widetilde{\sigma}} A_{\widetilde{\sigma}}^{-1}\beta)/\kappa}| \prod_{i\notin \widetilde{\sigma}\cup \{n+1\}} \left(
\sum_{k_i +m_i\in \NN} \frac{(C_2 |x^{(b_i)_A} \zeta^{ (b_i)_{n+1}}|)^{(k_i + m_i )}}{(k_i + m_i )!^{|b_i|} }\right)
\end{align*}
and for all $i\notin \widetilde{\sigma}\cup \{n+1\}$ we have:
\begin{equation}
\sum_{k_i +m_i\in \NN} \frac{(C_2 |x^{(b_i)_A} \zeta^{
(b_i)_{n+1}}|)^{(k_i + m_i )}}{(k_i + m_i )!^{|b_i|}
}\leq \left\{\begin{array}{ll} K_1 \operatorname{exp}(K_2 |C_2
x^{(b_{i})_A}|^{1/|b_i|}|\zeta|^{(b_i)_{n+1}/|b_i|}) & \mbox{ if }
(b_i)_{n+1}>0\\
\frac{1}{1-|C_2 x^{(b_i)_A} \zeta^{ (b_i)_{n+1}} |} & \mbox{
if } (b_i)_{n+1}\leq 0
\end{array} \right.
\end{equation} Take $U=\{ x \in \CC^n : |x^{(b_i)_A}|< R_i , \; i=1,\ldots ,n \}$ where $R_i>0$ can be chosen
arbitrarily large except when $(b_i)_{n+1}=|(b_i)_A|=0$ in which case we take $R_i<1/C_2$.  Then, since $(b_i)_{n+1}/|b_i|<1$ if $(b_i)_{n+1}>0$,
we have the result for for arbitrarily small $c_2 >0$ if we take $c_2>0$ and $R'>0$ big enough.\QED

\begin{proposition}\label{bounded-growth}
 $\varphi_B (x, \zeta)$ has an analytic continuation to an open set of the form $U\times S(\theta ,\delta)$,
where $U$ is
certain open set of $\CC^n$ and $S(\theta ,\delta)$ is a sector with bisecting direction $\theta$ and small
enough opening $\delta >0$.
Moreover, if $\Phi_{A_{\eta}}^{F,d-1}\subseteq \Phi_A^{F,d-1}$ then for arbitrarily small $c_2>0$ we can
chose $c_1 >0$ such that $|\varphi_B (x, \zeta)|\leq c_1 \operatorname{exp}(c_2 |\zeta |)$ for $(x,\zeta )\in U \times S(\theta ,\delta)$.
\end{proposition}

{\it Proof}\/.
To simplify the exposition we will first assume that $\alpha$ is
very generic.

Notice that $\eta \cup \{n+1\}$ is the set of (indices of) columns
of $A_B$ belonging to the hyperplane $H=\{ |A_{\sigma}^{-1}y|-
\frac{1}{|\operatorname{det}(A_{\sigma})|} y_{d+1} = 1 \}$ (see
Remark \ref{two-hyperplanes}) and we denote by $A'$ the submatrix
of $A_B$ consisting of these columns. Let $q$ be the cardinality of
$\eta$. Recall that $\varphi_B (x,\zeta)$ defines a holomorphic
function at each point of $U_{\sigma ,R}'$ (see Proposition
\ref{Gevrey-for-sigma-weight-vector} and (\ref{open-set})) and
notice  that $w_j = 0$ for all $j\notin \eta$.

We can write $$\varphi_B (x,\zeta)=\sum_{m\in \NN^{n-q} }
\varphi_{m} \dfrac{x_{\overline{\eta}}^m}{m!}$$ where
$\varphi_{m}=\varphi_m(x_\eta ,\zeta )$ is a holomorphic solution of
$H_{A'}(\beta - \sum_{i\notin \eta} m_i a_i,\alpha)$, which is
regular holonomic because all the columns of $A'$ belong to the
hyperplane $H$ \cite{hotta}. Let $W\subset \CC^{q+1}$ be the open
set such that $U_{\sigma ,R}'=W\times \CC^{n-q}$ (see
(\ref{open-set})), so that for all $m=(m_i)_{i\not\in \eta}\in
\NN^{n-q}$, $\varphi_m$ is holomorphic in $W$.

Take $Z=\{x_i =0: \; i\notin \eta \}$ and notice that we can identify $W$ with a
relative open subset of $Z$, i.e. with $W\times \{0\}=U_{\sigma ,R}'\cap Z$.

Recall that the singular locus of a hypergeometric system does not
depend on the  parameter but only on the matrix (see
\cite{Adolphson} and \cite{gkz1989}). In particular, since
$\varphi_m$ is convergent in $W$ for all $m$, we can consider the
analytic continuation of all the $\varphi_m$ along the same path
starting at a point in $W$ and avoiding the singular locus of the
hypergeometric system associated with $A'$. Let $c=(c_i )_{i\in
\eta} \in \CC^{q}$ be such that the complex line
$\{x_i = c_i :\; i \in \eta \}\cap Z$ (with coordinate $\zeta$)
intersects $W$ at nonsingular points of $H_{A'} (\beta ,\alpha)$.
Notice that this intersection is a relative open set in the complex
line. Let $\operatorname{Sing}(c)$ be the set of points $\zeta_0 \in
\CC\setminus \{0\}$ such that $(x_{\eta},\zeta )=(c ,\zeta_0)$ is a
singular point of $H_{A'}(\beta - \sum_{i\notin \eta} m_i a_i,\alpha)$.
$\operatorname{Sing}(c)$ is a finite set and thus $\Theta ( c)=\{
\arg (u ) :\; u \in \operatorname{Sing}(c) \}$ is also finite. As we
vary $c$ in a small open set $W'\subseteq \CC^{q}$, $\Theta
(c)$ is contained in a finite union of small intervals and we
can take $\theta$ such that for $\delta >0$ small enough, $(\theta -
\delta /2 ,\theta +\delta /2)\cap \Theta (x)=\emptyset$ for all
$x\in W'$. Hence we can consider the analytic continuation of each
$\varphi_m$ to an open set containing $W' \times S(\theta , \delta
)$.

We have extended $\varphi_B$ as a formal solution of $M_{A_B }(\beta ,\alpha )$ along $Z$, which
is convergent in some relative open set of $Z$. Thus, by the
constructibility of the solutions of a holonomic system in the sheaf
$\cO_{\widehat{X|Z}}/\cO_{X|Z}$ (see \cite{Mebkhout-positivite}) we have that the formal solution
constructed is convergent at any point of $W' \times S(\theta , \delta )\times \{\underline{0}\}$, thus
$\varphi_B (x,\zeta)$ can be analytically continued to an open set containing $W' \times S(\theta , \delta )
\times \{\underline{0}\}$.

Let us see that the analytic continuation of $\varphi_B (x,\zeta)$
satisfies a growth estimate near $\zeta =\infty$. Each $\varphi_{m}$
has polynomial growth since it is a solution of the regular
hypergeometric system $M_{A'}(\beta - \sum_{i\notin \eta} m_i
a_i,\alpha)$. Since $\alpha ,\beta$ are very generic, each
$\varphi_{m}$ can be written as a Nilsson series that converges in
certain open set (see e.g. \cite[Proposition 3.4.4]{SST}) which is a
linear combination of series of the form $\phi_{v(m)}(x_\eta ,\zeta
)$ (for some set of exponents $v(m)$ associated with simplices in
certain regular triangulation of the matrix $A'$) with support
$N_{v(m)}$ given by integer vectors in $\ker(A')$ with coordinates
sum equal to zero. The open set where the Nilsson series converge
depends on the regular triangulation of $A'$ that the simplices
belong to. We need to use a Nilsson series expression of
$\varphi_{m}$ that converges in points $(x_{\eta},\zeta)=(c ,\zeta)$
with $|\zeta | > R$ for a sufficiently large $R>0$. It is enough to consider a regular triangulation $T$ of $A_{\eta}$
and take $T'=\{\widetilde{\sigma}\cup \{n+1\}:\;  \widetilde{\sigma} \in T\}$ as the regular triangulation of
$A'$. By properties of regular triangulations, there is one regular triangulation $T$ of $A_{\eta}$ such that
there exists $c$ as above so that if $|\zeta | > R$ for a sufficiently large $R>0$ then $(c ,\zeta)$ belongs to
the domain
of convergence of the series $\phi_{\widetilde{v}}$ for any vector $\widetilde{v}$ associated with
$\widetilde{\sigma}\cup \{n+1 \}\in T'$.

The series expression of $\varphi_B$ via substitution of each
$\varphi_m$ by its Nilsson series expansion is a formal Nilsson
series (see e.g. \cite[Lemma 6.15]{Fer}). We know that this Nilsson series converges to $\varphi_B$ at points
in $W' \times S(\theta , \delta )\times \{\underline{0}\}$ close to
$\zeta =\infty$. By Lemma \ref{technical-lemma} it verifies the desired growth estimate.

Notice that the Nilsson series expansion of $\varphi_B$ near $\zeta =\infty$ also provides an analytic continuation
to points with $|\zeta |$ big enough and $x$ varying in certain open set of $\CC^n$ that contains $(c,0)$.

We have considered analytic continuations along paths contained in
$Z$. We can also extend $\varphi_B$ by analytic continuation along
paths from a point in $W \times \CC^{n-q}$ to a point near $\zeta
=\infty$ avoiding the singular locus of $M_{A_{B}}(\beta , \alpha
)$. If the starting point of the path is close to $(c,0)$ the
analytic continuation coincides with the Nilsson series close to
$\zeta =\infty$ and thus it satisfies the same growth estimate.

Finally, the parameter $\beta$ being very generic, the rank of $M_{A_B} (\beta_B )$ equals $vol(A_B)$ since
the set of exceptional parameters has codimension at least 2 \cite[Porism 9.5]{MMW}.
So, we can reduce the general case (when $\alpha$ is not
necessarily generic) to the previous one following the ideas of the
proof of \cite[Theorem 3.5.1]{SST} (see also the proof of \cite[Theorem 6.2.
]{Fer}).\QED

\begin{remark}\label{w-in-borel-revisited}
Using Proposition \ref{bounded-growth}, the results in Section
\ref{section:Borel} also hold for $w$ as in
(\ref{sigma-weight-vector}) instead of $w$ satisfying the assumption in
Theorem \ref{borel-thm} if we assume the additional condition
$\Phi_{A_{\eta}}^{F,d-1}\subseteq \Phi_A^{F,d-1}$ to hold. In particular, we
obtain an analogous version of Corollary \ref{corollary-asymptotic}. Let $\beta
\in \CC^n$ be very generic and let $\phi_v$ be a Gevrey solution of
$H_A (\beta )$ of index $s=1 +1/k>1$ with respect to some coordinate
subspace $Z\subseteq \CC^n$. Let $\sigma$ be the simplex which
$v$ is associated with, so we have that $\phi_v$ is also a
Gevrey series of index $s=1+1/k
>1$ with respect to $Y=\{x_i=0 :\; |A_{\sigma}^{-1} a_i |>1 \}
\supset Z$. Note that if we take $w$ associated with $\sigma$ as in
(\ref{sigma-weight-vector}) then $w (v+u)\in \NN$ for all $u\in
N_v$. By Proposition \ref{Gevrey-for-sigma-weight-vector} we have
that $t^{\alpha}\psi (x,t)=\phi_v (t^{w_1}x_1 ,\ldots, t^{w_n}x_n )$ is a
Gevrey series along $t=0$ of Gevrey index $s' = 1 + 1/\kappa$ with
$\kappa =|\det (A_{\sigma })|$. Let $\mathcal{S}[\psi](x,t)$ be the
$\kappa$-sum of $\psi (x ,t)$ with respect to $t$ in a direction $\theta \notin \Theta (x)$ for $x$
in certain open set $U$ small enough with compact closure. For any closed
subsector $\overline{S}$ of $S(\theta ,\alpha ,\rho)$ (see notations in Section \ref{section:Borel})
there exist constants $C>0, K>0$ such that the inequality
$|t^{\alpha}\mathcal{S}[\psi](c,t)-t^{\alpha} \psi_N ( c,t)| \leq C K^N \Gamma (1+N/\kappa ) |t|^N $ holds for
$t\in \overline{S}$, $c=(c_1,\ldots,c_n)\in U$ and any $N\in \NN$. Thus,
considering parametric curve $x(t)=(c_1 t^{w_1} , \ldots , c_n t^{w_n})$ then
$$\psi_N (x(t),1)=\sum_{u\in N_v , w
(v+u) \leq N - 1} \dfrac{[v]_{u_{-}}}{[v+u]_{u_{+}}} c^{v+u} t^{w (v+u)}
=t^{\alpha}\psi_N (c, t )$$ and $\mathcal{S}[\psi](x(t),1)= t^{\alpha}\mathcal{S}[\psi](c,t)$. Note that
$x$ tends to the point $x' \in Y$, with $x_i ' = x_i $ if
$|A_{\sigma}^{-1}a_i|\leq 1$, as $t$ tends to $0$.
\end{remark}

\begin{theorem}\label{theorem-asymptotic-expansion}

Let $\beta \in \CC^n$ be very generic and let $\phi_v$ be a Gevrey
solution of $H_A (\beta )$ of order $s=1 +1/k>1$ with respect
to a coordinate hyperplane $Y=\{x_i =0\}$. Let $\sigma$ be the simplex
which $v$ is associated with and take $w$ associated
with $\sigma$ as well. If $\Phi_{A_{\eta}}^{F,d-1}\subseteq \Phi_A^{F,d-1}$ then
for $\kappa =|\det (A_{\sigma })|$ we have that $\mathcal{S}[\psi ](x , 1)$
is a holomorphic solution of $M_A (\beta )$ and that for each $(x_1 ,\ldots ,x_{i-1},x_{i+1} ,\ldots, x_n)$ in
certain open set of $\CC^{n-1}$, $\phi_v (x)$ is a Gevrey asymptotic expansion of order $s$ of
$\mathcal{S}[\psi ](x, 1)$ with respect to $x_i=0$ in all but finitely many directions.
\end{theorem}

{\it Proof.-}
Assume for simplicity that the hyperplane is $Y=\{x_n =0\}$ and so $\sigma \subseteq \{1 ,\ldots ,n-1\}$.
We have that $w_i=0$ for $i=1 ,\ldots , n-1$ and $w_n =|\det (A_{\sigma })|(s-1)>0$ where $s=|A_{\sigma}^{-1} a_n |>1$ is
the Gevrey index of $\phi_v$ along $Y$.

By Remark \ref{w-in-borel-revisited} for $x(t)=(c_1 ,\ldots ,c_{n-1} , c_n t^{w_n})$ with
$t\in \overline{S}$ and $c\in U$, we have the inequality
$$|\mathcal{S}[\psi](x(t) ,1)- \psi_N (x(t),1)| \leq C
K^N \Gamma (1+N/\kappa ) |t|^N $$ for all $N\geq 0$.

Here we write $\psi_N (x,1)=\sum_{w_n m < N} f_m(x_1 ,\ldots , x_{n-1} ) x_n^m$. For integers of the form $N= w_n M$ with
$M\in \NN$ we have that
$t^N=(x_n(t) /c_n)^M$ and $N/\kappa = M (s-1)=M/k$ where $s=1 + 1/k$. Then we get the inequality
$$|\mathcal{S}[\psi](x(t) ,1)- \psi_N (x(t),1)| \leq C (K^{w_n}/|c_n|)^M \Gamma (1+ M/k) |x_n(t)|^M $$ for all $M\geq 0$.
This finishes the proof since we can assume by taken a smaller open set $U$ that $c_n \neq 0$ and that $|c_n|> C'$ for some constant $C'>0$.
\QED

\begin{example}
Put $A=(1 \; 2 \; 3)$, $\beta \in \CC$ and $w =(0,0,1)$. The vector
$v=(0,\beta /2 ,0)$ is an exponent of the $A$-hypergeometric system
$H_A (\beta )$ with respect to a perturbation of $w$ and so the
series
$$\psi (x,t)=\phi_v (x_1 , x_2 , t x_3 )=\sum_{m_1 , m_3 \geq 0, (m_1 + 3 m_3 )\in 2 \ZZ } \dfrac{[\beta /2]_{(m_1 + 3 m_3 )/2}}{m_1 ! m_3 !} x_1^{m_1 } x_2^{(\beta -m_1 - 3 m_3  )/2} x_3^{m_3} t^{m_3 }$$ is one of the series considered in the proof of
Theorem \ref{theorem-dimension-1} and it is a Gevrey solution of the
modified system $M_{A,w }(\beta )$ along $T$ with order $s=r+1=3/2$.
Notice that $s=r+1=3/2$ is the Gevrey index of $\psi (x,t)$ along
$T$ if and only if $\beta \notin 2 \NN$ (otherwise $\psi (x,t)$ is a
polynomial). Following Section \ref{section:Borel} but with our
vector $w$ (which does not satisfy the assumptions in Section
\ref{section:Borel} but is of the form \eqref{sigma-weight-vector}
for $\sigma = \{2\}$) we consider the Borel transform of $\psi$ with
index $\kappa = 1/r =2$:
$$\psi_B (x , \tau )=\sum_{m_1 , m_3 \geq 0, (m_1 + 3 m_3 )\in 2 \ZZ } \dfrac{[\beta /2]_{(m_1 + 3 m_3 )/2}}{m_1 ! m_3 ! \Gamma ( 1 + m_3 /2 )} x_1^{m_1 } x_2^{(\beta -m_1 - 3 m_3  )/2} x_3^{m_3} \tau^{m_3}.$$ This series defines a holomorphic function in
$\{(x,\tau)\in \CC^4 :\; | \frac{x_3 \tau}{x_2^{3/2}}|< \epsilon, x_1 , x_2 \neq
0 \}$ for $\epsilon>0$ small enough and it has an analytic continuation
with respect to $\tau$ to certain sector $S(\theta, \delta)$. Let us see that this analytic continuation has
polynomial growth in $\tau$. If $\varphi(x, z ):= \psi(x, t)|_{t = z^{1/2}}$ then its Borel
transform (with index $1$)
$$\varphi_B (x,\zeta)=\sum_{m_1 , m_3 \geq 0, (m_1 + 3 m_3 )\in 2 \ZZ } \dfrac{[\beta /2]_{(m_1 + 3 m_3 )/2}}{m_1 ! m_3 ! \Gamma ( 1 + m_3 /2 )} x_1^{m_1 } x_2^{(\beta -m_1 - 3 m_3  )/2} x_3^{m_3} \zeta^{m_3 / 2 }$$ is a solution of the
hypergeometric system associated  with $A_B$ and $(\beta ,0)$
defined on $\CC^4$ with coordinates $(x,\zeta)=(x_1 , x_2 , x_3 ,\zeta)$.
Notice that $\varphi_B$ has fractional powers in $\zeta$ but defines a
multivalued holomorphic function in $\{(x,\zeta )\in \CC^4 :\; |
\frac{x_3^2 \zeta}{x_2^3}|< \epsilon, \zeta \neq 0 \}$ for $\epsilon>0$ small enough.

It is clear that $\varphi_B (x,\zeta)$ is a linear combination of
series $\phi_{v} (x,\zeta )$ with $v\in \CC^{4}$ associated with the
simplex $\{2,4\}$ of $A_B$ (i.e., $A_B v =\beta_B$, $v_i\in \NN$ for
$i=1,3$). We have an analytic continuation of $\varphi_B (x,\zeta)$
to a point in the open set $\{(x,\zeta )\in \CC^4 :\; | \frac{x_3^2
\zeta}{x_2^3}|> R\}$ for $R>0$ big enough, which must be a linear
combination of series $\phi_{\widetilde{v}} (x,\zeta )$ with
$\widetilde{v}\in \CC^{4}$ associated with the simplex $\{3,4\}$
(this simplex alone determines a regular triangulation of the of
$A_B$ and the set of series $\phi_{\widetilde{v}}$ with
$\widetilde{v}$ associated with $\{3,4\}$ generates the space of
holomorphic solutions of $M_{A_B}(\beta_B )$ at any point in the
open set $\{(x,\zeta )\in \CC^4 :\; | \frac{x_3^2 \zeta}{x_2^3}|> R
\}$). We have that the columns of $B_{\{3,4\}}$ are $b_1
=(1,0,-1/3,-1/6), b_2 =(0,1,-2/3,-1/3)$ and we notice that $(b_1)_4
=-1/6 , (b_2)_4=-1/3 <0$. This implies that a series
$\phi_{\widetilde{v}} (x ,\zeta)$ with $\widetilde{v}$ associated
with $\{3,4\}$ has polynomial growth as $\zeta$ tends to infinity.
Thus, the analytic continuation of $\varphi_B (x,\zeta)$ close to
$\zeta =\infty$ also does.

As a consequence, the Borel sum  of $\psi$ with index $2$ given by
$$\mathcal{S}[\psi](x,t)=
\mathcal{L}^{\theta}_{2}\hat{\mathcal{B}}_{2}[\psi](x,t)=\int^{e^{i\theta}\cdot
\infty}_0 e^{-(\tau/t)^2} \psi_B (x , \tau )
d(\tau/t)^{2}$$ is a holomorphic solution of
$M_{A,w}(\widetilde{\beta })$ and $\mathcal{S}[\psi](x,1)$ is a
holomorphic solution of $M_A (\beta)$ which has an asymptotic
expansion $\psi (x,1)$ that is Gevrey of order $s=3/2$ along $x_3 =
0$.
\end{example}

\begin{example}
This example shows that the hypothesis in Proposition
\ref{bounded-growth} on the umbrella is necessary and that the bound
there is sharp. Take $$A=\left(\begin{array}{cccc}
           1 & 1 & 0 & \ell \\
           0 & 1 & 2 & 0\end{array}\right)$$
where $\ell>1$ is a rational number. We can consider $\ell\in \QQ$
by changing the lattice $\ZZ^2$ by the lattice $(\frac{1}{\ell}\ZZ)
\times \ZZ$. Then, the weight vector $w=(0,0,0,\ell-1)$ is
associated with $\sigma=\{1,2\}$ by the formula
(\ref{sigma-weight-vector}) and  $\kappa=|\det A_{\sigma}|=1$. Let
$\alpha\in \CC$ and $\beta \in \CC^2$ be very generic and let $v\in
\CC^4$ be a vector associated with $\sigma$ so that the series
$\phi_v(x)$ is a Gevrey solution of $M_A (\beta)$ along $x_4 =0$
with Gevrey index $s=1+r=\ell>1$. Then for $\psi_v (x,t) =
t^{-\alpha}\phi_v (x_1 ,x_2,x_3,t^{\ell-1} x_4)$ the Borel transform
$\varphi_B (x,\zeta )$ is convergent in the open set $\{(x ,\zeta ):
\; |\zeta^{\ell-1} |< \epsilon |x_1^\ell / x_4 |\}$ for some
$\epsilon>0$ small enough.

Moreover, it defines a holomorphic solution of $M_{A_B}(\beta_B )$
and then it is a linear combination of the set of series
$\phi_{v'}(x,\zeta)$ with $v'\in \CC^5$ associated with the simplex
$\{1,2,5\}$. Its analytic continuation to points in the open set $\{
(x ,\zeta ): \;|\zeta^{\ell-1} |>R |x_1^\ell / x_4 |\}$ for some
sufficiently large $R>0$ is a linear combination of the set of
series $\phi_{\widetilde{v}}$ with $\widetilde{v}\in \CC^5$
associated with the simplex $\{2,4,5\}$ of $A_B$. The column vectors
of $B_{\{2,4,5\}}$ are $b_1 =(1,0,0,-1/\ell, (1-\ell)/ \ell)$ and
$b_3=(0,-2,1,2/\ell,2(\ell-1)/\ell)$. Elements in the support
$N_{\widetilde{v}}$ of $\phi_{\widetilde{v}}$ are of the form $m_1
b_1 + m_3 b_3\in \ZZ^5$ with $m_1 , m_3 \in \NN$.   Thus any of
these series $\phi_{\widetilde{v}}(x,\zeta)$ is convergent in the
open set $\{ (x ,\zeta ): \;|\zeta^{\ell-1} |>R |x_1^\ell / x_4 |\}$
for some sufficiently large $R>0$ and since the last coordinate of
$b_3$ is $2(\ell-1)/\ell>0$ there is a subseries of
$\phi_{\widetilde{v}}$ (the one with monomials
$(x,\zeta)^{\widetilde{v}+m_3 b_3}$, $2m_3(\ell-1)/\ell \in \NN$)
such that the set of exponents of $\zeta$ in its monomials is
contained in $\widetilde{v}_5 +\NN$. The fact that $|b_3|=1>0$
guarantees that the coefficients of this subseries has the same type
of growth as $1/(m_1) !$ and thus the growth of this series, as
$\zeta$ tends to $\infty$ in certain sector, is equivalent to the
growth of $K \operatorname{exp}(C x_3 x_4^{2/\ell}
\zeta^{2(\ell-1)/\ell}/x_2^2)$ for some $K ,C>0$. Notice that for
$\ell>1$, we have that $\eta =\{1,2,4\}$ and the hypothesis
$\Phi_{A_{\eta}}^{F,d-1}\subseteq \Phi_A^{F,d-1}$ (required in Lemma
\ref{technical-lemma} and Proposition \ref{bounded-growth}) is
satisfied if and only if $1<\ell<2$. Thus, the bound
$|\phi_{\widetilde{v}}(x ,\zeta)|\leq c_1 \operatorname{exp}(c_2
|\zeta |)$ is satisfied for some $c_1 ,c_2 >0$ if and only if
$1<\ell \leq 2$ but for $\ell=2$ we cannot choose $c_2$ to be
arbitrarily small.
\end{example}


\begin{thebibliography}{99}

\bibitem{Adolphson} A.~Adolphson, \textit{A-hypergeometric functions and rings generated by
monomials}. Duke Mathematical Journal 73 (1994), n. 2, p. 269-290.

\bibitem{ACG-how-1996} A.~Assi, F.J.~Castro-Jim\'enez  and  J.-M.~Granger, \textit{How to
calculate the slopes of a D- module.} Compositio Math., 104 (1996)
107-123.

\bibitem{B}
W.~Balser, From divergent power series to analytic functions,
Lecture Notes in Mathematics {\bf{1582}},
Springer, 1994.

\bibitem{bfp} C.~Berkesch, J.~Forsg\r{a}rd, M.Passare,
Euler-Mellin Integrals and $A$-Hypergeometric Functions,
arxiv:1103.6273v2.

\bibitem{Bernstein70} J.~Bernstein, \newblock{\em Modules over the ring of differential operators.
Study of the fundamental solutions of equations with constant
coefficients.} \newblock {Funkcional. Anal. i Prilo\v zen.  5
(1971), no. 2, 1--16.}

\bibitem{beukers} F.~Beukers, Monodromy of $A$-Hypergeometric Functions,
arxiv:1101.0493v2.

\bibitem{castro-takayama} F.J.~Castro-Jim\'{e}nez, N.~Takayama, \newblock{\em Singularities of the hypergeometric D-module associated with a monomial curve}, Trans. Am. Math. Soc. 355(9) (2003) 3761-3775.

\bibitem{DMM} A.~Dickenstein, F.~Mart\'{\i}nez and L.~Matusevich, \textit{Nilsson solutions for irregular A-hypergeometric systems}. Rev. Mat. Iberoam. 28 (3) (2012) 723-758.

\bibitem{esterov-takeuchi} A.~Esterov, K.~Takeuchi,
{\em Confluent A-hypergeometric functions and rapid decay homology cycles,}
arxiv:1107.0402, 2011.

\bibitem{Fer} M.~C.~Fern\'andez Fern\'andez, \emph{Irregular hypergeometric
$D$-modules}, Adv. Math. 224 (2010) 1735-1764.

\bibitem{fer-cas-1} M.C.~Fern\'{a}ndez-Fern\'{a}ndez, F.J.~Castro-Jim\'{e}nez, {\em Gevrey solutions of irregular hypergeometric systems in two variables}, J. Algebra 339 (2011) 320-335.

\bibitem{fer-cas-2} M.C.~Fern\'{a}ndez-Fern\'{a}ndez, F.J.~Castro-Jim\'{e}nez, {\em Gevrey solutions of the irregular hypergeometric system associated with an affine monomial curve}, Trans. Am. Math. Soc. 363(2) (2011) 923-948.

\bibitem{ggz1987}
I.~M.~Gel'fand, M.~I.~Graev and A.~V.~Zelevinsky, Holonomic systems
of equations and series of hypergeometric type, Dokl. Akad. Nauk
SSSR 295 (1987), no. 1, 14--19.

\bibitem{gkz1989} I.~M.~Gel'fand, A.~V.~Zelevinsky, M.~M.~Kapranov,
Hypergeometric functions and toral manifolds.
Functional Analysis and its Applications {\bf 23} (1989), 94--106.

\bibitem{gkz-advance} I.M.~Gel'fand, M.~Kapranov, A.~Zelevinsky,
Generalized Euler Integrals and $A$-Hypergeometric Functions,
Advances in Mathematics {\bf 84} (1990), 255--271.

\bibitem{hartillo1} M.I.~Hartillo, {\em Hypergeometric slopes of codimension 1}, Rev. Mat. Iberoam. 19 (2) (2003) 455-466.

\bibitem{hartillo2} M.I.~Hartillo, {\em Irregular hypergeometric systems associated with a singular monomial curve},  Trans. Amer. Math. Soc. 357 (2005), no. 11, 4633\^{a}4646.

\bibitem{hotta} R.~Hotta, Equivariant D-modules, arxiv:math/9805021, 1998.

\bibitem{Hra85} J.~Hrabowski, {\em Multiple hypergeometric functions and simple Lie groups $SL$ and $Sp$}. SIAM J. Math. Anal. {\bf 16}, (1985), 876--886.

\bibitem{Laumon} G.~Laumon, \textit{Transformations canoniques et sp\'ecialisation pour les D-modules filtr\'es}, Ast\'erisque  {\bf 130}, (1985), 56--129.

\bibitem{Laurent} Y.~Laurent, \textit{Polyg\^one de Newton et b-fonctions pour les modules microdiffer\'entiels}, Ann. Sci. \'Ecole Norm. Sup. (4) {\bf 20} (1987), 391--441.

\bibitem{Laurent-Mebkhout-pa-pa} Y.~Laurent and  Z.~Mebkhout, \textit{Pentes
alg\'ebriques et pentes analytiques d'un $\cD$-module}. Annales
Scientifiques de L'E.N.S. $4^e$ s\'erie, tome 32, n. 1 (1999) 39-69.

\bibitem{MMW} L.~F.~Matusevich, E.~Miller and U.~Walther, \textit{Homological
methods for hypergeometric families}. J. Amer. Math. Soc. 18 (2005),
4, 919-941.

\bibitem{Mebkhout-positivite} Z.~Mebkhout, \textit{Le th\'eor\`eme de positivit\'{e} de l'irr\'egularit\'{e} pour les D-modules}. Progress in. Math., vol. 88, 83-132. Birkh\"{a}user, Boston, 1990.

\bibitem{ohara-takayama} K.~Ohara and N.~Takayama,  \textit{Holonomic rank of
A-hypergeometric differential-difference equations}. J. Pure Appl.
Algebra 213 (2009), no. 8, 1536-1544.

\bibitem{SST}
M.~Saito, B.~Sturmfels, and N.~Takayama, \emph{Gr\"obner
{D}eformations of {H}ypergeometric {D}ifferential {E}quations},
Springer--Verlag, Berlin, 2000.

\bibitem{SW} M.~Schulze and U.~Walther, \textit{Irregularity of
hypergeometric systems via slopes along coordinate subspaces}. Duke
Math. Journal 142, {\bf 3} (2008), 465-509.

\bibitem{SW2} M.~Schulze and U.~Walther, \textit{Irregularity of
hypergeometric systems via slopes along coordinate subspaces}. arXiv:math/0608668v3 [math.AG].

\bibitem{Sturmfels} B.~Sturmfels, Gr\"obner Bases and Convex Polytopes. American Mathematical Society,
Univ. Lectures Series, Band 8. Providence, Rhode Island, 1996.

\bibitem{takayama2009}
N.~Takayama, Modified ${\mathcal A}$-Hypergeometric Systems, Kyushu
Journal of Mathematics {\bf 63} (2009),  113--122.

\end{thebibliography}
\end{document}